\documentclass[11pt]{article}
\usepackage{amsmath}
\usepackage{amsmath,amssymb,amsfonts,amsthm,fancyhdr}
\usepackage{epsfig,graphicx,picins,picinpar,subfigure}
\usepackage{pstricks}
\usepackage{fancyvrb}
\usepackage[numbers,sort&compress]{natbib}
\begin{document}

%%%%%%%%%%%%%%%%%%%%%%%%%%%%%%%%%%%%%%%%%%%%%%%%%%%%%%%%%%%%%%%%
% 标题,作者, 通信地址定义
%%%%%%%%%%%%%%%%%%%%%%%%%%%%%%%%%%%%%%%%%%%%%%%%%%%%%%%%%%%%%%%%
%%\begin{CJK}{GBK}{song}
\title{On the Infinitude of Some Special Kinds of Primes \\------ {\small\it Dedicated to the memory of my mother}}
\author{Shaohua Zhang}
\date{{\small  School of Mathematics, Shandong University,
Jinan, China, 250100 \\E-mail address:
shaohuazhang@mail.sdu.edu.cn}}
 \maketitle
%%%%%%%%%%%%%%%%%%%%%%%%%%%%%%%%%%%%%%%%%%%%%%%%%%%%%%%%%%%%%%%%
%  英文文摘要
%%%%%%%%%%%%%%%%%%%%%%%%%%%%%%%%%%%%%%%%%%%%%%%%%%%%%%%%%%%%%%%%
%\begin{center}\textbf{Abstract}\quad \end{center}
%\vspace{-1em}
%\begin{abstract}
\textbf{Abstract:}   The aim of this paper is to try to establish a
generic model for the problem that several multivariable
number-theoretic functions represent simultaneously primes for
infinitely many integral points. More concretely, we introduced
briefly the research background-the history and current
situation-from Euclid's second theorem to Green-Tao theorem. We
analyzed some equivalent necessary conditions that irreducible
univariable polynomials with integral coefficients represent
infinitely many primes, found new necessary conditions which perhaps
imply that there are only finitely many Fermat primes,  obtained an
analogy of the Chinese Remainder Theorem, generalized Euler's
function, the prime-counting function and Schinzel-Sierpinski's
Conjecture and so on.  Nevertheless, this is only a beginning and it
miles to go. We hope that number theorists consider further it.


\textbf{2000 MR  Subject Classification:}\quad 11A41, 11A99
%\end{abstract}
%\vspace{3mm}
%%%%%%%%%%%%%%%%%%%%%%%%%%%%%%%%%%%%%%%%%%%%%%%%%%%%%%%%%%%%%%%%
%  正文由此开始
%%%%%%%%%%%%%%%%%%%%%%%%%%%%%%%%%%%%%%%%%%%%%%%%%%%%%%%%%%%%%%%%

\section{Research background---the history and current
situation---from Euclid's second theorem to Green-Tao theorem}
\setcounter{section}{1}\setcounter{equation}{0}
%%%%%%%%%%%%%%%%%%%%%%%%%%%%%%%%%%%%%%%%%%%%%%%%%%%%%%%%%%%%%%%%
From ancient to modern times, the study of the infinitude of some
special kinds of primes has been one of the most important topics in
Number Theory. People usually ask the following questions:

\vspace{3mm}1, Are there infinitely many Fermat primes?  Fermat
primes are primes of the form $2^{2^x}+1$.

2, Are there infinitely many Mersenne primes? Mersenne primes are
primes of the form $2^x-1$, where $x$ is also a prime.

3, Are there infinitely many twin primes?

4, Are there infinitely many primes of the form $x^2+1$?

5, Are there infinitely many Sophie Germain primes?  A prime $p$ is
called a Sophie Germain prime if $2p+1$  is also prime.

And so on.

\vspace{3mm}Mathematicians throughout history have been fascinated
by these problems. However, they are still unanswered.  Euclid [1]
proved firstly the following result.

\vspace{3mm}\noindent {\bf Euclid's second theorem:~~}%
There are infinitely many primes.

\vspace{3mm}Anyone who likes Number Theory must like Euclid's second
theorem. In his book \emph{The book of prime number records }[2],
Paulo Ribenboim cited nine and a half proofs of Euclid's second
theorem. In this paper, we listed the references of fifteen distinct
proofs again, see [3--17].

\vspace{3mm}Clearly, using Euclid's method, the ancient Greeks can
also prove that there are infinitely many primes of the form $4k-1$
or $6k-1$. Using properties of quadratic residues, it is easy to
prove that there are infinitely many primes of the form $4k+1$ or
$6k+1$. Cyclotomic polynomials [18] can be used to prove that there
are infinitely many primes of the form $ak+1$. In 2004, Yoo, Jisang
[19] gave another elementary proof of the infinitude of primes of
the form $ak+1$. Especially, in 2005, Robbins Neville [20] gave a
simple proof of the infinitude of primes of the form $3k+1$.

\vspace{3mm}Naturally, a more general problem on primes in
arithmetic progressions seems that there should be infinitely many
primes of the form $a+bn$, where $a$ and $b$ are integers satisfying
$(a,b)=1$, and either $a\neq0, b>0$, or $a=0, b=1$. After the time
of Euclid, there have been no great improvements on this problem in
about 2000 years. Until 1837, using L-series and analytic methods,
Dirichlet [21] solved thoroughly it.

\vspace{3mm}\noindent {\bf Dirichlet's theorem:~~}%
There are infinitely many primes of the form $a+bn$, where $a$ and
$b$ are integers satisfying $(a,b)=1$, and either $a\neq0, b>0$, or
$a=0, b=1$.

\vspace{3mm}This is a classical and most important theorem which is
perceived as a milestone of the study on the infinitude of some
special kinds of primes. In the 1890's, de la Vall\'{e}e Poussin
[22] showed further that the number of such primes not exceeding a
large number $x$ is asymptotic to $x/{\varphi(b)\log x}$ as
$x\rightarrow\infty$, where $\varphi(.)$ is Euler's function.

\vspace{3mm} Clearly, the question of existence of infinitely many
primes in arithmetic progressions can be regard as the question of
existence of infinitely many prime values of linear polynomials. In
1857, Bouniakowsky [23] considered the case of nonlinear polynomials
and stated a conjecture below.

\vspace{3mm}\noindent {\bf Bouniakowsky's conjecture:~~}%
If $f(x)$ is an irreducible polynomial with integral coefficients,
positive leading coefficient and degree at least 2, and there does
not exist any integer $n>1$ dividing all the values $f(k)$ for every
integer $k$, then $f(x)$ is prime for an infinite number of integers
$x$.

\vspace{3mm}Concerning the simultaneous values of several linear
polynomials, Dickson [24] stated the following conjecture in 1904:

\vspace{3mm}\noindent {\bf Dickson's conjecture:~~}%
Let $s\geq1$, $f_i(x)=a_i+b_ix$ with $a_i$ and $b_i$ integers,
$b_i\geq1$ (for $i=1,...,s$ ). If there does not exist any integer
$n>1$ dividing all the products $\prod_{i=1}^{i=s}f_i(k)$, for every
integer $k$, then there exist infinitely many natural numbers $m$
such that all numbers $f_1(m),...,f_s(m)$ are primes.

\vspace{3mm}In 1958, by studying the consequences of Bouniakowsky's
conjecture and Dickson's conjecture, A. Schinzel and W. Sierpinski
[25] got the following conjecture:

\vspace{3mm}\noindent {\bf Schinzel-Sierpinski conjecture (H hypothesis):~~}%
Let  $s\geq1$, and let $f_1(x),...,f_s(x)$ be irreducible
polynomials with integral coefficients and positive leading
coefficient. If there does not exist any integer $n>1$ dividing all
the products $\prod_{i=1}^{i=s}f_i(k)$, for every integer $k$, then
there exist infinitely many natural numbers $m$ such that all
numbers $f_1(m),...,f_s(m)$ are primes.

\vspace{3mm}For some details on primes represented by univariate
polynomials, see also [26--32]. As for the case of primes
represented by polynomials in few variables, it is very complicated
and precise conjectures do not seem to have been formulated in the
literature for multivariable polynomials still less univariable
number-theoretic functions or multivariable number-theoretic
functions. However, some notable results on the question of
existence of infinitely many prime values of bivariate polynomials
have been obtained by using sieve methods.

\vspace{3mm}The problem goes back to Fermat who proved that there
are infinitely many primes of the form $x^2+y^2$. E. Schering [33]
and H. Weber [34] proved that every primitive binary quadratic form
(positive if definite) with discriminant different from a perfect
square represents infinitely many primes. In 1969, Motohashi, Yoichi
[35] proved that there are infinitely many primes of type
$x^2+y^2+1$. In the early 1970's, as an improvement of results of
Bredihin B. M., Linnik Ju. V. and Motohashi Yoichi [35-37], H.
Iwaniec [38, 39] obtained the significant asymptotic formula of the
number of primes represented by a primitive quadratic polynomial. In
1997, Fouvry, Etienne and Iwaniec, Henryk [40] proved that there are
infinitely many primes of type $x^2+y^2$, where $x$ is a prime
number.

\vspace{3mm}In the above-mentioned sequences of polynomial values in
which it has been proved there are infinitely many primes, there are
$\gg x/(\log x)^c $ elements of the sequence up to $x$, for some
fixed $c>0$. Below are two great results on some bivariate
polynomials can take on infinitely many prime values.

\vspace{3mm}In 1998, Friedlander, John and Iwaniec, Henryk [41]
proved that $x^2+y^4$ takes on a prime value for $\sim z^{3/4}/\log
z$ values $\leq z $, which implies that there are infinitely many
primes of the type $x^2+y^4$ . It is a "monumental
breakthrough"---reviewed by Andrew Granville.

\vspace{3mm}In 2001, Heath-Brown, D. R. [42] proved that $x^3+2y^3$
takes on a prime value for $\sim z^{2/3}/\log z$ values $\leq z $,
which implies that there are infinitely many primes of the type
$x^3+2y^3$. It is "one of the major landmarks of analytic number
theory"--- reviewed by G. Greaves.

\vspace{3mm}After the work of Friedlander, John and Iwaniec, Henryk
in 1998 and Heath-Brown, D. R. in 2001, maybe, the next goal of this
line of research is to prove that Landau's first conjecture [66] is
true. Namely, there are infinitely many primes of the form $x^2+1$.
It should be interesting to see.

\vspace{3mm}As for the primes of other forms, such as Wilson primes,
Wieferich primes, regular primes, NSW-primes, the primes of form
$\frac{10^x-1}{9}$ and so on, see many papers or books, for example
[43-62].

\vspace{3mm}Finally, we introduce the famous work of Ben Green,
Terence Tao and Tamar Ziegler to close this section. In 2004, Ben
Green and Terence Tao [63] proved the following brilliant result:

\vspace{3mm}\noindent {\bf Green-Tao theorem:~~}%
The sequence of prime numbers contains arbitrarily long arithmetic
progressions.

\vspace{3mm} Green-Tao theorem is a great support to Dickson's
conjecture and this deep and important result has brought a very
significant impact in studying primes. "It is a landmark
contribution to additive number theory."---reviewed by Tamar
Ziegler. Recently, they further gave important consideration and
profound analysis on Dickson's conjecture [64]. In 2006, Terence Tao
and Tamar Ziegler [65] extended Green-Tao theorem to polynomial
progressions via the Bergelson-Leibman polynomial Szemer\'{e}di
theorem.

\vspace{3mm}Based on the aforementioned rich achievements and
advancements, and also due to the fact that the universe has been
governed by the same laws, we believe that it is possible to
establish a generic model for the problem that several multivariable
number-theoretic functions represent simultaneously primes for
infinitely many integral points. It will be the main aim of this
paper. Nevertheless, this is a very intractable task. It seems that
the author can not finish well. Our work is only a beginning. We
hope that number theorists consider further it.

\vspace{3mm}Next, let's begin with the simplest case that an
irreducible univariable polynomial with  integral coefficients
represents infinitely many primes.

%%%%%%%%%%%%%%%%%%%%%%%%%%%%%%%%%%%%%%%%%%%%%%%%%%%%%%%%%%%%%%%%
\section{Necessary conditions that an irreducible univariable
polynomial with integral coefficients represents infinitely many
primes}
%%%%%%%%%%%%%%%%%%%%%%%%%%%%%%%%%%%%%%%%%%%%%%%%%%%%%%%%%%%%%%%%
In this paper, we always restrict that a $k$-variables
number-theoretic function $f(x_1,...,x_k)$ is a map from $N^k$ to
$Z$. Moreover, we assume that $f(x_1,...,x_k)$ is a continuous
function on $R^k$, where $R$ is the set of all real numbers.
Specially, an irreducible univariable polynomial $f(x)$ is a map
from $N$ to $Z$. Of course, a prime number is  positive. We do not
consider negative primes.

\vspace{3mm}Let $f(x)$ be a univariable polynomial with  integral
coefficients, we further assume
 that $f(x)$ is not a constant.
 Note that pairwise distinct primes are pairwise relatively
 prime. Thus, we get a natural necessary condition that $f(x)$ represents infinitely many primes.

\vspace{3mm}\noindent {\bf Necessary condition A:~~}%
There exists an infinite sequence of positive integers
$x_1,x_2,...,x_k,...$ such
 that$f(x_1),f(x_2),...,f(x_k),...$ are pairwise relatively prime, moreover
 $f(x_1)>1,f(x_2)>1,...,f(x_k)>1,...$.

\vspace{3mm}\noindent {\bf Proposition 1:~~}%
Necessary condition A implies the following necessary conditions B,
C and D. Moreover, B, C and D are equivalent.

\vspace{3mm}\noindent {\bf Necessary condition B:~~}%
For any positive integer $m>1$, there exists a positive integer
 $x$ such that $\gcd (f(x),m)=1$.

\vspace{3mm}\noindent {\bf Necessary condition C:~~}%
For any positive integer $m>1$, there exists a positive integer
 $x$ such that $m$ does not divide $f(x)$. Namely, there does not exist a positive integer $m>1$,
 such that for any positive integer $x$, $m$ divides $f(x)$.

\vspace{3mm}\noindent {\bf Necessary condition D:~~}%
For any prime $p$, there exists a positive integer
 $x$ such that $\gcd (f(x),p)=1$. Namely, there does not exist a prime $p$, such that
 for any positive integer $x$, $p$ divides $f(x)$.

\vspace{3mm}\noindent {\bf Proof of Proposition 1:~~}%
Clearly, $A\Longrightarrow B\Longrightarrow C\Longrightarrow D$.
Next, we prove that $D\Longrightarrow B$. For any positive integer
$m>1$, we write $m=\prod_{i=1}^{i=k} {p_i}^{e_i}$. By $D$, there
exists a positive integer
 $a_i$ such that $\gcd (f(a_i),{p_i}^{e_i})=1 $  for $1\leq i\leq k$.  By Chinese Remainder
 Theorem, there
exists a positive integer $x$ such that $x\equiv a_i (\mod
{p_i}^{e_i})$. Note that $f(x)$ is a  polynomial with integral
coefficients. ( Here, $f(x)$ need not be irreducible. ) Hence,
$f(x)\equiv f(a_i) (\mod {p_i}^{e_i})$ and $\gcd (f(x),m)=1$.

\vspace{3mm}Obviously, if $f(x)$ represents infinitely many primes,
then we must have the following:

\vspace{3mm}\noindent {\bf Necessary condition E:~~}%
For any positive integer $m>1$, there exists a positive integer
 $x$ such that $\gcd (f(x),m)=1$ and $f(x)>1$.

\vspace{3mm}\noindent {\bf Necessary condition F:~~}%
For any positive integer $m>1$, there exists a positive integer
 $x$ such that $f(x)>1$ and $m$ does not divide $f(x)$.

\vspace{3mm}\noindent {\bf Necessary condition G:~~}%
For any prime $p$, there exists a positive integer
 $x$ such that $\gcd (f(x),p)=1$ and $f(x)>1$.

\vspace{3mm} \noindent {\bf Proposition 2:~~}%
Necessary conditions A and E  are equivalent, however, they and F
(resp.
 G) are not always equivalent.

\vspace{3mm}\noindent {\bf Proof of Proposition 2:~~}%
Let $f(x)$ be a  polynomial with integral coefficients. If the
leading coefficient of $f(x)$ is positive, then A, E, F and G are
equivalent by the idea of proof in Proposition 1. Now we consider
the case that the leading coefficient of $f(x)$ is negative.
Clearly, in this case, we still have: $A\Longrightarrow E$. Next, we
prove that $E\Longrightarrow A$. Since the leading coefficient of
$f(x)$ is negative, hence $f(x)$ at most represents finitely many
positive integers. Denote the product of these positive integers by
$M$. By Necessary conditions E, let $m=2$, there exists a positive
integer
 $x$ such that $f(x)>1$ and $\gcd (f(x),2)=1$ . This implies that
 $f(x)$ can
represent a positive integer greater than 1. Therefore, $M>1$. By E
again, there exists a positive integer
 $x$ such that $f(x)>1$ and $\gcd (f(x),M)=1$. But $M$ is the product of
all positive integers which can be represented by $f(x)$. It is
impossible. So, A and E  are equivalent. By considering
$f(x)=-x^2+6$, it is easy to prove that $f(x)=-x^2+6$ implies F and
G, but A and F (resp.
 G) are not always equivalent because $f(x)=-x^2+6$ can not represent
 infinitely many positive integers and $f(x)=-x^2+6$ does not imply A.

\vspace{3mm} Based on Proposition 2, Necessary condition A and E
will become our main
 interest in future study.

\vspace{3mm}\noindent {\bf Corollary 1:~~}%
Let $f(x)$ be a polynomial with integral coefficients, then A and
the following condition are equivalent: the leading coefficient of
$f(x)$ is positive, and there does not exist any integer $n>1$
dividing all the values $f(k)$ for every integer $k$.

\vspace{3mm}\noindent {\bf Remark 1:~~}%
Some people call Necessary condition C (resp. D) Bunyakovsky's
property.

\vspace{3mm}\noindent {\bf Remark 2:~~}%
Our work in this section shows that there are several equivalent
forms of Bouniakowsky's conjecture. For instance, if $f(x)$ is an
irreducible polynomial with integral coefficients, and for any
positive integer $m>1$, there exists a positive integer
 $a$ such that $\gcd (f(a),m)=1$ and $f(a)>1$, then $f(x)$ represents
 infinitely many primes.

\vspace{3mm}\noindent {\bf Remark 3:~~}%
Generalizing our work to the generic cases, one could obtain several
equivalent forms of Dickson's conjecture even Schinzel-Sierpinski's
Conjecture. For example, let  $s\geq1$, and let $f_1(x),...,f_s(x)$
be irreducible polynomials with integral coefficients, if there
exists an infinite sequence of positive integers
$x_1,x_2,...,x_k,...$ such
 that $\prod_{i=1}^{i=s}f_i(x_1),...,\prod_{i=1}^{i=s}f_i(x_k),...$ are
 pairwise relatively prime, moreover
 $f_i(x_1)>1,f_i(x_2)>1,...,f_i(x_k)>1,...$ for $i=1,...,s$, then there exist
infinitely many natural numbers $m$ such that all numbers
$f_1(m),...,f_s(m)$ are primes.

\vspace{3mm}We do not know Dickson, Schinzel and Sierpinski whether
noticed these equivalent forms. It seems that they focused on
Bunyakovsky's property and believed that if a univariable polynomial
$f(x)$ with integral coefficients and the positive leading
coefficient has Bunyakovsky's property, then $f(x)$ represents
infinitely many primes. Namely, for any univariable polynomial
$f(x)$ with  integral coefficients and the positive leading
coefficient, Bunyakovsky's property of $f(x)$ is the sufficient and
necessary condition that $f(x)$ represents infinitely many primes.

\vspace{3mm}Unfortunately, these conjectures are open for many
years. It is time to reconsider them. On one hand, one maybe ask: is
Bunyakovsky's property of $f(x)$ enough to determine that $f(x)$
represents infinitely many primes? On the other hand, how to
generalize Schinzel-Sierpinski's Conjecture to the cases of
multivariable polynomials with  integral coefficients even
multivariable number-theoretic functions?

\vspace{3mm}Let $f_1(x_1,...,x_k),...,f_s(x_1,...,x_k)$ be $s$
multivariable number-theoretic functions from $N^k$ to $Z$. Assuming
that $f_1(x_1,...,x_k),...,f_s(x_1,...,x_k)$  represent
simultaneously primes for infinitely many integral points
$(x_1,...,x_k)$. Now we generalize Necessary condition A to the
generic case as follows.

\vspace{3mm}\noindent {\bf Necessary condition H:~~}%
There exists an infinite sequence of integral points
$(x_{11},...,x_{k1})$, ..., $(x_{1i},...,x_{ki})$, ... such that
$\prod_{j=1}^{j=s}f_j(x_{11},...,x_{k1})$,...,
$\prod_{j=1}^{j=s}f_j(x_{11},...,x_{ki})$,... are pairwise
relatively prime and $f_j(x_{11},...,x_{ki})>1$ for each $i$ and
$j$.

\vspace{3mm}Similarly, Necessary condition H and the following
necessary condition I  are equivalent:

\vspace{3mm}\noindent {\bf Necessary condition I:~~}%
For any positive integer $m>1$, there exists an integral point
 $X=(x_1,...,x_k)$ such that $\gcd (\prod_{i=1}^{i=s}f_i(X),m)=1$ and $f_i(X)>1$ for $1\leq i\leq s$.

\vspace{3mm}As we aforementioned, Necessary condition H should be
viewed as a natural necessary condition. Based on this observation,
we believe that there is always a common necessary condition that
any multivariable number-theoretic functions represent
simultaneously primes for infinitely many integral points. Surely,
at least, it is not weaker than the natural necessary condition and
can be called the maximum necessary condition. Once adding
appropriate conditions, it perhaps leads to the sufficient condition
that multivariable number-theoretic functions represent
simultaneously primes for infinitely many integral points.
Therefore, it is possible to generalize Schinzel-Sierpinski's
Conjecture.

\vspace{3mm}We also find that using the natural necessary condition
H (resp. I) is more convenient than using Bunyakovsky's property
when we treat the multivariable cases, in which we have not the
definition of leading coefficient even irreduciblity. Moreover, our
work will show that the natural necessary condition perhaps is the
maximum necessary condition when
$f_1(x_1,...,x_k),...,f_s(x_1,...,x_k)$ are multivariable
polynomials with  integral coefficients. For details, see next
several sections.

%%%%%%%%%%%%%%%%%%%%%%%%%%%%%%%%%%%%%%%%%%%%%%%%%%%%%%%%%%%%%%%%
\section{Find new necessary conditions}
%%%%%%%%%%%%%%%%%%%%%%%%%%%%%%%%%%%%%%%%%%%%%%%%%%%%%%%%%%%%%%%%
{\bf Why do we need to find new necessary
conditions?~~}%
Note that the number-theoretic function $2^{2^x}+1$ implies the
natural necessary condition A. Numbers of the form $2^{2^x}+1$ are
called Fermat numbers. Primes of the form $2^{2^x}+1$ are Fermat
primes. Eisenstein proposed as a problem in 1844 the proof that
there are infinite number of Fermat primes [2]. Nevertheless, Hardy
and Wright [67] conjectured that the number of Fermat primes is
finite, although they did not give any reasons and explanations. By
factoring Fermat number, many people believe that the conjecture in
[67] holds. So far, people do not find a new Fermat primes except
for the first four Fermat primes as follows: 5, 17, 257, 65537. If
let $x=0$, then 3 is viewed as a Fermat prime. But we restricted
that a number-theoretic function is a map from $N$ to $Z$ in Section
2. Therefore, here, we do not consider 3.

\vspace{3mm}Historically, the problem that the number-theoretic
function $f(x)=2^{2^x}+1$  represents primes were first studied by
Pierre de Fermat, who conjectured that $f(x)=2^{2^x}+1$ are prime
for all $x\in N\cup \{0\}$. Unfortunately, in 1732, his conjecture
was refuted by Leonhard Euler. Euler showed that $f(5)=4294967297
=641\times 6700417$. Euler proved that every prime factor of $f(x)$
must have the form $k\times 2^{x+1}+1$. For $x=5$, this means that
the only possible factors are of the form $64k+1$. Euler found the
factor 641 when $k=10$. Lucas refined Euler's result: Any prime
divisor of $f(x)$ is of the form $k\times 2^{x+2}+1$ whenever $x>1$.

\vspace{3mm}According to R. P. Brent [92]: "The complete
factorization of Fermat numbers $f(6)$, $f(7)$, ..., has been a
challenge since Euler's time. Because the $f(x)$ grow rapidly in
size, a method which factors $f(x)$ may be inadequate for $f(x+1)$,
No Fermat primes larger than $f(4)$ are known, and a probabilistic
argument makes it plausible that only a finite number of $f(x)$
(perhaps only 3, 5, 17, 257, 65537) are prime."  As of 2008 it is
known that $f(x)$ is composite for $5 \leq x \leq 32$, although
complete factorizations of $f(x)$ are known only for $0 \leq x \leq
11$, Below is the list of complete factorizations:

\vspace{3mm}$f(6)=18446744073709551617=274177\times 67280421310721$
[88].

$f(7)=59649589127497217\times p22$, where $p22$ is a prime which has
22 decimal digits [89].

$f(8)=1238926361552897\times p62$ [90].

$f(9) = 2424833\times p49\times p99$  [91].

$f(10)=45592577\times 6487031809\times p40\times p252$ [92].

$f(11)=319489\times 974849\times p21\times p22\times p564$ [92].

\vspace{3mm}Thus, if the conjecture in [67] holds, then the natural
necessary condition A is too weak to make us know more information
on the infinitude of some special kinds of primes, and it should be
strengthened. Here, we always assume that there is the maximum
necessary condition that any multivariable number-theoretic
functions represent simultaneously primes for infinitely many
integral points.

\vspace{3mm}Another reason, although it seems reluctant, we give it
as follows: Considering the infinitude, the technical definition of
limit occurred to us. If $f_1(x_1,...,x_k),...,f_s(x_1,...,x_k)$
represent simultaneously primes for infinitely many integral points
$(x_1,...,x_k)$, then there should exist a constant $c$ such that
for every positive integer $m>c$, there exists an integral point
$(x_1,...,x_k)$ such that $f_1(x_1,...,x_k),...,f_s(x_1,...,x_k)$
are coprime to $m$ simultaneously. Based on some heuristic
observations, for example, by refining Necessary condition E as
follows: for a sufficiently large constant $c$ and for any positive
integer $m>c$, there exists a positive integer
 $x$ such that $\gcd (f(x),m)=1$ and $m>f(x)>1$, we would like to restrict the values of
$f_1(x_1,...,x_k)$, ... , $f_s(x_1,...,x_k)$ in
$$Z_m^* = \{x\in N | 1\leq x<m, \gcd(x,m)=1 \}$$ in order to know more information that
$f_1(x_1,...,x_k),...,f_s(x_1,...,x_k)$ maybe take on infinitely
many prime values.

\vspace{3mm}The third reason, let's look back Necessary condition E
again. Necessary condition E states that for an integral polynomial
$f(x)$ and for any positive integer $m>1$, there exists a positive
integer
 $x$ such that $\gcd (f(x),m)=1$ and $f(x)>1$. Therefore,
 if for such an integral polynomial $f(x)$, Necessary condition E satisfies,
 then there must exist the least positive integer
 $n$ such that $\gcd (f(n),m)=1$ and $f(n)>1$. Denote this least positive
 integer $n$ by $S_f(m)$, then for a sufficiently large constant $c$ and for any positive integer
$m>c$, $S_f(m)<(\frac{m}{L(f)})^{\frac{1}{d}}$, where $d$ is the
degree of $f(x)$, and, $L(f)$ is the leading coefficient of $f(x)$.

\vspace{3mm}For given polynomial $f(x)$ and every sufficiently large
$m$, Estimating the upper bound of $S_f(m)$ is an interesting
question. For example, let $f(x)=x$, then $S_f(m)<m^{\frac{1}{2}}$
when $m>30$ by Bonse's inequalities [68-69]. Moreover, this result
can be refined as follows: for any given positive integer $k$, there
is a constant $c_k$, when $m>c_k$, $S_f(m)<m^{\frac{1}{k}}$. As
another example, generalizing $f(x)$ to the case of number-theoretic
functions and defining similarly $S_f(m)$, let $f(x)=2^x-1$, then
$S_f(m)<\log_2m$ when $m>21$ [70]. In fact, when $f(x)=2^x-1$, the
meaning of $S_f(m)$ is definite. when $\gcd(m,2)=1$,
$\gcd(m,2^{\varphi (m)+1}-1)=1$; when $\gcd(m,2)=2$, we write $m=2^e
t$ with $\gcd(t,2)=1$, then $\gcd(m,2^{\varphi (t)+1}-1)=1$. When
$f(x)=2^{2^x}+1$, for any positive integer $m$, we have
$\gcd(m,2^{2^m}+1)=1$ because any prime divisor $p$ of $2^{2^m}+1$
is of the form $k2^{m + 2} + 1$ whenever $m$ is greater than one. So
$S_f(m)\leq m$ when $f(x)=2^{2^x}+1$.

\vspace{3mm}$S_f(m)$ also can be generalized to the generic case:
Let $f_1(x_1,...,x_k)$, ... ,$f_s(x_1,...,x_k)$ be multivariable
number-theoretic functions. If for any positive integer $m$, there
exists an integral point $(y_1,...,y_k)$ such that
$f_1(y_1,...,y_k)>1,...,f_s(y_1,...,y_k)>1$ are all coprime to $m$,
then there must exist the shortest integral vector $X=(x_1,...,x_k)$
such that $f_1(X)>1,...,f_s(X)>1$ are all coprime to $m$. Denote
this shortest integral vector by $S_{f_1,...,f_s}(m)$. Then,
$S_{f_1,...,f_s}(m)$ is the generalization of $S_f(m)$.

\vspace{3mm}Estimating the upper bound of $S_f(m)$ also leads to
strengthen Necessary condition E as the aforementioned.

\vspace{3mm}Certainly, making the decision of strengthening
Necessary condition E should always take a risk. We must verify the
sequences of functional values in which it has been proved there are
infinitely many primes implies that there exist a constant  $c$ such
that for every positive integer $m>c$, there exist an integral point
such that those corresponding functional values are all in $Z_m^* $.

\vspace{3mm}To begin with, noticed that if $a>1$ is the smallest
integer such that $ \gcd(a,m)=1$, then  $a$ is a prime when $m>2$.
Namely, there exists a constant  $c=2$ such that for every positive
integer $m>c$, there is a prime in  $Z_m^* $. Thus, we proved the
case of $s=k=1$ with $f(x)=x$.

\vspace{3mm}In additionally, note that $\pi_{a,b,x}\sim
\frac{x}{\varphi(b)\log x}$ as $x\rightarrow\infty$, where
$\pi_{a,b,x}$ is the number of prime of the form $a+bx$ with $b>0,
\gcd (a,b)=1$. This implies that there is a positive constant  $c$,
when  $m>c$, we have $\pi_{a,b,x}>1+\log_2 m$. But $m$ has at most
$[\log_2 m]$ distinct prime factors. Hence, there is always a prime
of the form $a+bx$ in $Z_m^* $ and we proved the case $s=k=1$ with
$f(x)=ax+b$.

\vspace{3mm}Last but not the least, using the similar method, one
can show respectively that there is a positive constant  $c$, when
$m>c$, the number of prime of the form $f(x,y)=x^2+y^2+1$,
$f(x,y)=x^3+2y^3$, $f(x,y)=x^2+y^4$ and so on $>1+\log_2 m$ in the
cases of $s=1$ and $k=2$. It follows immediately the desired
consequences. Combining with the above discussions, one could
conjecture the following:

\vspace{3mm}\noindent {\bf Conjecture 1:~~}%
Let $f_1(x_1,...,x_k),...,f_s(x_1,...,x_k)$  be $s$ multivariable
number-theoretic functions. If
$f_1(x_1,...,x_k),...,f_s(x_1,...,x_k)$ represent simultaneously
primes for infinitely many integral points $(x_1,...,x_k)$, then
there is always a constant $c$  such that for every positive integer
$m>c$, there exists an integral point $(y_1,...,y_k)$ such that
$f_1(y_1,...,y_k)>1,...,f_s(y_1,...,y_k)>1$ are all in $Z_m^* $.

\vspace{3mm}\noindent {\bf Remark 4:~~}%
The conjecture is only a necessary condition not a sufficient
condition. For example, by Bonse's inequalities [68-69], one can
prove that every positive integer $m>30$, there are positive
integers of the form $x^2$ in $Z_m^* $. But, $x^2$ never represents
a prime.

\vspace{3mm}Now, we prove that Conjecture 1 implies that there are
only finitely many  Fermat primes. In fact, if there are infinitely
many Fermat primes, then by Conjecture 1, there is always a constant
$c$ such that for every positive integer  $m>c$, there exists a
positive integer $n$ such that $f(n)=2^{2^n}+1$ in $Z_m^* $. Let
$m=\prod_{i=0}^{i=k-1}f(i)=\prod_{i=0}^{i=k-1}(2^{2^i}+1)> c$.
Clearly, there is always such a positive integer $k$ because $c$ is
a constant. Hence, we must have $n\geq k$ when
$m=\prod_{i=0}^{i=k-1}(2^{2^i}+1)$ and $2^{2^n}+1$ in $Z_m^* $. Thus
$f(n)\geq f(k)$. Note that $f(k)=m+2$. Therefore, we have $f(n)\geq
m+2$. But, it is impossible since $2^{2^n}+1$ in $Z_m^* $ implies
that $2^{2^n}+1<m$ by the meaning of notation $Z_m^* $.

\vspace{3mm}In like manner, this follows immediately Conjecture 1
which implies also that there are only finitely many primes of the
form $n^n+1$. Based on the same reason, we maybe foresee that there
are only finitely many prime values for several iterative functions.
For example, maybe, there are only finitely many prime numbers in
the sequence: $p_1=2^2-1=3, p_2=2^3-1=7,  p_3=2^7-1=127,
p_4=2^{127}-1, p_5=2^{p_4}-1,....$

\vspace{3mm}A clear sense is of that such a sequence is so sparse
that it can not guarantee that there is always a constant $c$ such
that for every positive integer $m>c$, there exists a positive
integer  $n$ such that $p_n$ is in $Z_{p_1p_2...p_{n-1}}^* $.

\vspace{3mm}Besides, we demand Conjecture 1 to test that the
infinitude of some special kinds of primes such as Twins primes,
safe primes (co-Sophie-Germain primes), Mersenne primes and so on,
which are markedly infinitely many. In [70], we proved that the
several number-theoretic functions ($f(x)=x, g(x)=x+2$; $f(x)=x,
g(x)=2x+1$; $f(x)=x, g(x)=2^x-1$) which perhaps represent
simultaneously infinitely many primes imply Conjecture 1. Moreover,
by the following quantitative form of Schinzel-Sierpinski's
Conjecture---Bateman-Horn's conjecture [71], if $f_1(x),...,f_s(x)$
are polynomials with integral coefficients, and represent
simultaneously infinitely many primes, then Conjecture 1 holds.

\vspace{3mm}\noindent {\bf Bateman-Horn's conjecture:~~}%
Let  $s\geq1$, and let $f_1(x),...,f_s(x)$ be irreducible
polynomials with integral coefficients and positive leading
coefficient. If there does not exist any integer $n>1$ dividing all
the products $\prod_{i=1}^{i=s}f_i(k)$, for every integer $k$, and
for every integer $m>1$, the number $Q(m)$ of integers $1\leq n\leq
m$ such that $f_1(n),...,f_s(n)$ are all primes is about
$$C_{f_1,...,f_s}\frac{1}{\prod_{i=1}^{i=s}d_i}\sum_{n=2}^{n=m}\frac{1}{(\log
n)^s}\sim C\frac{1}{\prod_{i=1}^{i=s}d_i}\int_2^m\frac{dt}{(\log
t)^s},$$ where $d_i =\deg f_i(x)$,
$C=C_{f_1,...,f_s}=\prod_p\frac{1-\omega (p)/p}{(1-1/p)^s}$ is a
very complicated constant and $\omega (p)$ is the number of
solutions $x$, $0\leq x\leq p-1$, of the congruence
$f_1(x)...f_s(x)\equiv 0 (\mod p)$.

\vspace{3mm}Roughly speaking, the number $Q(m)$ is about
$C\frac{1}{\prod_{i=1}^{i=s}d_i}\frac{m}{(\log m)^s}$, which of
course, implies Conjecture 1 when $f_1(x),...,f_s(x)$ are
polynomials with integral coefficients. Namely, if polynomials
$f_1(x),...,f_s(x)$ with integral coefficients represent
simultaneously primes for infinitely many integers $x$, then there
is always a constant $c$ such that for every positive integer $m>c$,
there exists an integers $y$ such that $f_1(y)>1,...,f_s(y)>1$ are
all in $Z_m^* $.

\vspace{3mm}\noindent {\bf Remark 5:~~}%
Friedlander John and Granville Andrew [79-81] showed that
Bateman-Horn's asymptotic formula does not always hold and there are
infinitely many different polynomials of given degree which take
either significantly more or significantly less prime values than
expected. However, we believe that Conjecture 1 holds without a
proviso when $f_1(x),...,f_s(x)$ in Conjecture 1 are irreducible
polynomials with integral coefficients and positive leading
coefficient. In our another paper \emph{Notes on Dickson's
Conjecture}, we have proved strictly that Conjecture 1 holds when
$f_1(x),...,f_s(x)$ in Conjecture 1 are all linear polynomials with
integral coefficients and positive leading coefficient, for the
details, see [101]. Furthermore, in [101], we  generalize Dickson's
Conjecture to the multivariable case or a system of affine-linear
forms on $N^k$.  In [102], we give Dickson's conjecture on $Z^n$ and
obtain an equivalent form of Green-Tao's conjecture [64].

\vspace{3mm}Anyway, like the $\varepsilon-\delta$ definition of
limit, Conjecture 1 maybe provides us with another mathematical
description for the infinitude of some special kinds of primes.

\vspace{3mm}Conjecture 1 leads to the generalizations of Euler's
function and the prime-counting function, see Section 4. It also
yields an analogy of Chinese Remainder Theorem, see Section 5. But,
Conjecture 1 is based on the finiteness of  Fermat primes which is
unproved yet. Everyone is unwilling to see its unreliable basis. So,
we only hope that one keeps it in his mind. Conjecture 1 maybe will
lead to some correct conjectures.

\section{Generalizations of Euler's function and the prime-counting
function} Euler's totient function $\varphi (n)$ is a very important
number-theoretic function and defined to be the number of positive
integers $x$  less than $n$ which are relatively prime to $n$.
$\varphi (n)= \# \{x\in N|\gcd (x,n)=1, x<n\}=n\prod_{p|n}(1-1/p)$.
If we look upon $x$ as the value of number-theoretic function
$f(x)$, then when $f(x)=x$, we have $\varphi (n)= \# \{f(x)\in
N,\gcd (f(x),n)=1, f(x)<n |x\in N\}$. Thus, let $f(x)$ be a
number-theoretic function, then one can generalize Euler's totient
function as follows: $\Phi_f (n)= \# \{f(x)\in Z_n^* |x\in N\}$.

\vspace{3mm}More generally, let
$f_1(x_1,...,x_k),...,f_s(x_1,...,x_k)$ be $s$ multivariable
number-theoretic functions from $N^k$ to $Z$. one can generalize
further $\Phi_f (n)$ as follows: $$\Phi_{f_1,...,f_s} (n)= \#
\{f_1(X)\in Z_n^* ,...,f_s(X)\in Z_n^*|X=(x_1,...,x_k)\in N^k\}.$$

Now, we generalize another important number-theoretic function
--- the prime-counting function $\pi (x)$, which is the number of primes less
than or equal to some real number $x$. Note that pairwise distinct
primes are pairwise relatively prime. Consider the number-theoretic
function $f(x)=x$. For any given positive integer $x>1$, consider a
special sub-set $H$ of $\{1,2,...,x\}$ as following: $\forall a\in
H$, we have $a>1$, and $\forall a\neq b \in H$, we also have $\gcd
(a, b)=1$. Namely, the elements of $H$ are pairwise relatively
prime.

\vspace{3mm}Denote the set of all such sub-sets of $\{1,2,...,x\}$
by $M$. Thus, $M=\{H\subseteq \{1,2,...,x\}|\forall a\neq b \in H,
\gcd (a,b)=1,\forall a\in H, a>1\}$. Clearly, $\pi
(x)=\max_{H\subseteq M} \{\#H \}$. Namely, $\pi (x)$ can be viewed
as the largest among the cardinality of all sub-sets (in which each
element exceeds 1 and pairwise distinct elements are pairwise
relatively prime) of $\{1,2,...,x\}$.

\vspace{3mm}Now, let $f(x)$ be a generic number-theoretic function.
Let $H$ be any sub-set of the image of $f$. Consider the set
$M=\{H\subseteq \{1,2,...,x\}|\forall f(a)\in H, f(a)>1,\forall
f(a)\neq f(b) \in H, \gcd (f(a),f(b))=1\}$. Let
$\Pi_f(x)=\max_{H\subseteq M} \{\#H \}$. Then, $\Pi_f(x)$ can be
viewed as the generalization of $\pi (x)$. Denote the number of
distinct prime factors of $x$ by $\omega (x)$. If we have
$\Pi_f(m)>\omega (m)$, then there is a positive integer $a $ such
that $f(a)$ is in $Z_m^*$, and $\Phi_f (m)\geq 1$.

\vspace{3mm}More generally, let
$f_1(x_1,...,x_k),...,f_s(x_1,...,x_k)$ be $s$ multivariable
number-theoretic functions, consider the set $M=\{H\subseteq
\{1,2,...,x\}|\forall f(X)\in H, f(X)>1,\forall f_1(X)\neq f_1(Y)
\in H,..., f_s(X)\neq f_s(Y) \in H, \gcd (\prod
_{i=1}^{i=s}f_i(X),\prod _{i=1}^{i=s}f_i(Y))=1\}$, where integral
points $X,Y$ should be viewed as vectors.

\vspace{3mm}Let $\Pi_{f_1,...,f_s}(x)=\max_{H\subseteq M}\{\#H \}$.
Then, $\Pi_{f_1,...,f_s}(x)$ can be viewed as the generalization of
$\pi (x)$. Thus, if $f_1(x_1,...,x_k),...,f_s(x_1,...,x_k)$
represent simultaneously primes for infinitely many integral points
$X=(x_1,...,x_k)$, then, we must have
$\Pi_{f_1,...,f_s}(x)\rightarrow \infty$ as $x\rightarrow \infty$.

\vspace{3mm}Similarly, if $\Pi_{f_1,...,f_s}(m)>\omega (m)$ then
$\Phi_{f_1,...,f_s} (m)\geq 1$.

\vspace{3mm}Using sieve theory [72-75], one could obtain some
asymptotic formulae of $\Phi_{f_1,...,f_s} (m)$ and
$\Pi_{f_1,...,f_s} (m)$. This should become the subject of future
publications.

\section{An analogy of Chinese Remainder Theorem}
Chinese Remainder Theorem [76] states that for given a system of
simultaneous linear congruences  $x\equiv a_i (\mod n_i)$ for
$i=1,2,...,k$ and for which $n_i$ are pairwise relatively prime
positive integers, where $a_i$ are integers, then this linear system
has a unique solution modulo $n=\prod_{i=1}^{i=k}n_i$. Particularly,
for $i=1,2,...,k$, if $\gcd (a_i,n_i)=1$, then, this linear system
has a unique solution $x$ in $Z_n^*$.

\vspace{3mm}Chinese Remainder Theorem is the greatest theorem in
ancient China in my eyes. And it is a very theorem which was named
after a unique nation. It is one of the jewels of mathematics and
contains in a third-century AD book The Mathematical Classic by Sun
Zi by Chinese mathematician Sun Tzu. It reflects a perfect
combination of beauty and utility. The famous Fast Fourier Transform
can be even viewed as a special case of its. That is because Chinese
Remainder Theorem can be generalized over generic rings and Fourier
Transform formula $f\rightarrow (f(\omega ^0),...,f(\omega
^{n-1}))(f'\rightarrow \frac{1}{n}(f'~(1),...,f'(\omega ^{-n+1})))$
is exactly viewed as the isomorphism $C[x]/{x^n-1} \simeq
C[x]/{x-\omega}\times ...\times C[x]/{x-\omega ^n}$ which is
essentially Chinese Remainder Theorem. People said that "it is
difficult to image what would happen if there was no Fast Fourier
Transform in modern communications". This will enable us to learn
better the significance of Chinese Remainder Theorem. In this
section, we will give an analogy of its.

\vspace{3mm}Let us look back the proof of Theorem 2 in [70]. In
order to prove that there is always a constant $c$, such that when
$n>c$, there exists $x\in Z_n^*$ and $2x+1\in Z_n^*$ with $x>1$, our
method is to prove firstly that there exists $x\in Z_a^*$ and
$2x+1\in Z_a^*$, to prove secondly that there exists $y\in Z_b^*$
and $2y+1\in Z_b^*$ with $\gcd (a,b)=1$, to prove lastly that there
exists $z\in Z_{ab}^*$ and $2z+1\in Z_{ab}^*$. This is exactly
viewed as Chinese Remainder Theorem which states essentially that if
there is an integer in $Z_a^*$, and there is an integer in $Z_b^*$
with $\gcd (a,b)=1$, there is an integer in $Z_{ab}^*$. We hope
certainly that this can be generalized to generic cases as follows.

\vspace{3mm}\noindent {\bf An analogy of Chinese Remainder Theorem:~~}%
Let $f_1(x_1,...,x_k)$, ... , $f_s(x_1,...,x_k)$ be multivariable
polynomials with  integral coefficients. If
$f_1(x_1,...,x_k),...,f_s(x_1,...,x_k)$  represent simultaneously
primes for infinitely many integral points,  and if $\gcd (a,b)=1$
and there exist integral point $(x_1,...,x_k)$ and $(y_1,...,y_k)$
such that $f_1(x_1,...,x_k)>1,...,f_s(x_1,...,x_k)>1$ are all in
$Z_a^*$, and $f_1(y_1,...,y_k)>1,...,f_s(y_1,...,y_k)>1$ are all in
$Z_b^*$, then there exists an integral point $(z_1,...,z_k)$ such
that $f_1(z_1,...,z_k)>1,...,f_s(z_1,...,z_k)>1$  are all in
$Z_{ab}^*$.

\vspace{3mm}Here, we must explain why the condition that
"$f_1(x_1,...,x_k),...,f_s(x_1,...,x_k)$  represent simultaneously
primes for infinitely many integral points" is necessary. That is
because if the number of primes is finite, then Chinese Remainder
Theorem is false [77], namely, Chinese Remainder Theorem implies
Euclid's second theorem. In fact, $f(x)=x^3+1$ has not this property
because it does not represent infinitely many primes. For example,
$f(1)=2\in Z_9^*$ and $f(2)=9\in Z_{10}^*$, but there is not a
positive integer $x$ such that $f(x)=x^3+1\in Z_{90}^*$. Thus, we
also obtain another necessary condition that
$f_1(x_1,...,x_k),...,f_s(x_1,...,x_k)$ represent simultaneously
primes for infinitely many integral points.

\vspace{3mm}We also find that $f(n)=2^x-1$ satisfies this necessary
condition [70]. However, $f(n)=2^{2^n}+1$ does not satisfy this
necessary condition. For example, $5\in Z_{51}^*$ and $17\in
Z_{5\times 257}^*$. But, there is not a Fermat number in
$Z_{51\times 5 \times 257}^*$. Does it imply possibly that there are
only finitely many Fermat primes again? The answer perhaps is no, at
least we have a reason, due to the fact that we only consider the
case that $f_1(x_1,...,x_k),...,f_s(x_1,...,x_k)$ are multivariable
polynomials with  integral coefficients. We do not consider the case
of generic number-theoretic functions. In fact, in the case of
generic number-theoretic functions, we have not such an analogy of
Chinese Remainder Theorem. For instance, let $f(x)=\left\{
\begin{array}{c}
2,1\leq x\leq 2 \\
3,3\leq x\leq 39 \\
\lbrack x/3],x\geq 40\\
\end{array}
\right. $. Clearly, $f(x)$ represents infinitely many primes. But it
has not the similar property of Chinese Remainder Theorem. For
instance, $f(1)=2\in Z_3^*$ and $f(3)=3\in Z_4^*$, but there is not
a positive integer $x$ such that $f(x)\in Z_{12}^*$. (Is there a
counterexample of the analogy of Chinese Remainder Theorem when
$f(x_1,...,x_k)$ is a continuous function on $R^k$?) By this
example, one maybe ask: is Conjecture 1 true? We do not assert the
answer now. But it is possible to lead to generalize
Schinzel-Sierpinski's Conjecture.

\vspace{3mm}\noindent {\bf Remark 6:~~}%
Conjecture 1 and the analogy of Chinese Remainder Theorem should be
equivalent when $f_1(x_1,...,x_k),...,f_s(x_1,...,x_k)$ in
Conjecture 1 are multivariable polynomials with  integral
coefficients.

\vspace{3mm}\noindent {\bf Remark 7:~~}%
Note that if there are primes in $Z_a^*$ and $Z_b^*$ respectively,
then there is primes in $Z_{ab}^*$ when $\gcd (a,b)=1$. Similarly,
if multivariable polynomials $f_1(x_1,...,x_k)$, ... ,
$f_s(x_1,...,x_k)$ with  integral coefficients represent
simultaneously primes for infinitely many integral points, and if
$\gcd (a,b)=1$ and there are integral points $(x_1,...,x_k)$ and
$(y_1,...,y_k)$ such that $f_1(x_1,...,x_k),...,f_s(x_1,...,x_k)$ in
$Z_a^*$ are all primes, and $f_1(y_1,...,y_k)$, ... ,
$f_s(y_1,...,y_k)$ in $Z_b^*$ are all primes, then maybe, there
exists an integral point $(z_1,...,z_k)$ such that
$f_1(z_1,...,z_k),...,f_s(z_1,...,z_k)$  in $Z_{ab}^*$ are all
primes. This is a very interesting problem on the analogy of Chinese
Remainder Theorem, and in the simple case, we have: if $\gcd
(m,n)=1$, $\gcd (a,b)=1$ with $b>1$ and $a+bx\in Z_m^*$is prime, and
$a+by\in Z_n^*$ also is prime, then there exists a prime of the form
$a+bz$ in $Z_{mn}^*$.

\vspace{3mm}\noindent {\bf Remark 8:~~}%
When the paper is written here, we feel that it is not difficult to
generalize Schinzel-Sierpinski's Conjecture to the case of
multivariable polynomials with  integral coefficients. It will be
not a pure speculation anymore and become a somewhat reasonable
conjecture. For details, see Section 6.

\section{Generalizing Schinzel-Sierpinski's Conjecture
to the case of multivariable polynomials} A possible generalization
of Schinzel-Sierpinski's Conjecture is the following:

\vspace{3mm}Let $f_1(x_1,...,x_k),...,f_s(x_1,...,x_k)$ be
multivariable polynomials with  integral coefficients, if
$f_1(x_1,...,x_k),...,f_s(x_1,...,x_k)$ are irreducible over
$Q[x_1,...,x_k]$, and there is always a constant $c$  such that for
every positive integer $m>c$, there exists an integral point
$(y_1,...,y_k)$ such that $f_1(y_1,...,y_k)>1$, ... ,
$f_s(y_1,...,y_k)>1$ are all in $Z_m^* $, then
$f_1(x_1,...,x_k),...,f_s(x_1,...,x_k)$  represent simultaneously
primes for infinitely many integral points $(x_1,...,x_k)$.

\vspace{3mm}However, we do not do this. On one hand, there are many
puzzles on the factorization in $Q[x_1,...,x_k]$, and maybe, the
word "irreducible" can not explain more. On the other hand, in order
to generalize it to the more generic case such as number-theoretic
functions, we need a dependable condition to replace the
"irreducible" condition.

\vspace{3mm}For this goal, let us look back on the work of M. Ram
Murty [78]: let  $f(x)=\sum_{i=0}^{i=m}a_ix^i$ be a polynomial of
degree $m$ in $Z[x]$ and set $H=\max_{0\leq i\leq m-1}|a_i/a_m|$, If
$f(n)$ is prime for some integer $n\geq H+2$, then $f(x)$ is
irreducible in $Z[x]$. Based on the work of M. Ram Murty and our
aforehand analysis, we give the following conjecture.

\vspace{3mm}\noindent{\bf  Conjecture 2:~~}%
Let $f_1(x_1,...,x_k),...,f_s(x_1,...,x_k)$ be multivariable
polynomials with  integral coefficients, if there is a positive
integer $c$ such that for every positive integer $m\geq c$, there
exists an integral point $(y_1,...,y_k)$ such that
$f_1(y_1,...,y_k)>1,...,f_s(y_1,...,y_k)>1$ are all in $Z_m^* $, and
there exists an integral point $(z_1,...,z_k)$ such that
$f_1(z_1,...,z_k)\geq c,...,f_s(z_1,...,z_k)\geq c$ are all primes,
then $f_1(x_1,...,x_k),...,f_s(x_1,...,x_k)$  represent
simultaneously primes for infinitely many integral points
$(x_1,...,x_k)$.

\vspace{3mm}\noindent {\bf Remark 9:~~}%
Let $f_1(x_1,...,x_k),...,f_s(x_1,...,x_k)$ be multivariable
polynomials with  integral coefficients from $N^k$ to $Z$, then the
following  conditions are equivalent:

\vspace{3mm}$(H)$: If there exists an infinite sequence of integral
points $(x_{11},...,x_{k1})$, ... , $(x_{1i},...,x_{ki}),...$ such
that
$\prod_{j=1}^{j=s}f_j(x_{11},...,x_{k1}),...,\prod_{j=1}^{j=s}f_j(x_{11},...,x_{ki}),...$
are pairwise relatively prime and $f_j(x_{11},...,x_{ki})>1$ for
each $i$ and $j$.

\vspace{3mm}$(J)$: If there is a positive integer $c$ such that for
every positive integer $m\geq c$, there exists an integral point
$(y_1,...,y_k)$ such that
$f_1(y_1,...,y_k)>1,...,f_s(y_1,...,y_k)>1$ are all in $Z_m^* $.

\vspace{3mm}Thus, we deduce an equivalent form of
Schinzel-Sierpinski conjecture: Let $s\geq1$, and let
$f_1(x),...,f_s(x)$ be irreducible polynomials with integral
coefficients, if there is a positive integer $c$ such that for every
positive integer $m\geq c$, there exists a positive integer $a$ such
that $f_1(a)>1,...,f_s(a)>1$ are all in $Z_m^* $, then there exist
infinitely many natural numbers $m$ such that all numbers
$f_1(m),...,f_s(m)$ are primes.

\vspace{3mm}Conjecture 2 should view as the sufficient and necessary
condition that multivariable polynomials
$f_1(x_1,...,x_k),...,f_s(x_1,...,x_k)$  with  integral coefficients
 represent infinitely many primes.

\vspace{3mm}\noindent {\bf Remark 10:~~}%
Can Conjecture 2 be generalized similarly to the generic case of
number-theoretic functions? The author would like to keep it in mind
because this problem is unattackable now.  For example, let
$f_1=3\times 2^{x-1}-1,f_2=3\times 2^x-1, f_3=9\times 2^{2x-1}-1$,
then, do $f_1,f_2, f_3$ represent simultaneously primes for
infinitely many $x$? Another example, do $g_1=8x+5,g_2=x^3+2,
g_3=2^x-1$ represent simultaneously primes for infinitely many $x$?
Particularly, does $h(x)=2^x+x$ represent primes for infinitely many
$x$? And so on.

\vspace{3mm}In the author's eyes, it perhaps is easy to give a
sufficient condition that multivariable number-theoretic functions
$f_1(x_1,...,x_k),...,f_s(x_1,...,x_k)$ represent infinitely many
primes, but it is difficult to give its sufficient and necessary
condition. On this problem, we will try to present a plausible
proposal in Section 9.

\vspace{3mm}\noindent {\bf Remark 11:~~}%
Conjecture 2 leads to the following significative problem:

\vspace{3mm}Let $f_1(x_1,...,x_k),...,f_s(x_1,...,x_k)$ be
number-theoretic functions. Assume that there is a positive integer
$c$ such that for every positive integer $m\geq c$, there exists an
integral point $(y_1,...,y_k)$ such that
$f_1(y_1,...,y_k)>1,...,f_s(y_1,...,y_k)>1$ are all in $Z_m^* $, and
there exists an integral point $(z_1,...,z_k)$ such that
$f_1(z_1,...,z_k)\geq c,...,f_s(z_1,...,z_k)\geq c$  are all primes.
Since $f_1(x_1,...,x_k),...,f_s(x_1,...,x_k)$ represent
simultaneously primes for some integral point $(x_1,...,x_k)$, hence
we can denote the least prime represented simultaneously by
$P_{f_1,...,f_s}$. A significative problem is to estimate the upper
bound of $P_{f_1,...,f_s}$.

\vspace{3mm}Historically, this problem is one of important topics in
Number Theory. In the simplest case, denote $p(l,k)$ the least prime
in the arithmetic progression $l+kn$ with $(l,k)=1$, where $n$ runs
through the positive integers, and let $p(k)$ be the maximum value
of $p(l,k)$  for all $l$ satisfying $(l,k)=1$ and $1\leq l\leq k$.
Linnik proved that there exist positive $C$ and $L$ such that
$p(k)<Ck^L$. Heath-Brown proved that $p(k)<Ck^{5.5}$ [83]. On the
problem of the least prime in an arithmetic progression, Chinese
mathematicians and Chengdong Pan, Jingrun Chen, Jianmin Liu and Wei
Wang et al. made great contributions, see [93-100]. In the case of
irreducible polynomials with degree $>1$, McCurley Kevin S., Adleman
Leonard M., Odlyzko Andrew M. [27, 82, 84] obtained important
results.

\vspace{3mm}\noindent {\bf Remark 12:~~}%
Conjecture 2 is the first to mention the existence of primes among
the conjectures which conjecture the infinitude of some special
kinds of primes. Of course, if one wants to prove that the
infinitude, firstly, he must prove the existence. Unfortunately, it
is a critical difficulty. In next section, we go on with this
problem.

%%%%%%%%%%%%%%%%%%%%%%%%%%%%%%%%%%%%%%%%%%%%%%%%%%%%%%%%%%%%%%%%
\section{The existence of some special
kinds of primes}
%%%%%%%%%%%%%%%%%%%%%%%%%%%%%%%%%%%%%%%%%%%%%%%%%%%%%%%%%%%%%%%%
We begin with Euclid in this section. In his beautiful proof of the
infinitude of primes, Euclid must know the existence of primes. Of
course, the existence of primes is very clear. So he omitted the
proof of the existence of primes and supposed that there are only
finitely many primes, say $k$ of them, which denoted by
$p_1,...,p_k$ and constructed directly the number $1+
\prod_{i=1}^{i=k}p_i$ which leads to the contradiction.

\vspace{3mm}As we know, it is very difficult to prove the existence
of some special kinds of primes. For example, for every $k\geq 1$,
we even do not know whether there are always primes $p$ and $q$ such
that $p-q=2k$ or not. Namely, we do not know whether $f(x)=x$ and
$g(x)=x+2k$ represent simultaneously primes for some integer $x$ and
every $k$ so far.

\vspace{3mm}If one does not know whether there are some special
kinds of primes or not, can he prove their infinitude? This problem
goes back to Euler. By Euler's identity $\sum_{n=1}^\infty n^{-s}=
\prod_p (1-p^{-s})^{-1}$, one could prove that $\sum_p
p^{-s}\rightarrow \infty$ as $s\rightarrow 1$ which implies the
existence of infinitely many primes. Based on Euler's idea,
Dirichlet introduced Characters and proved further $\sum_{p\equiv
a(\mod b)} p^{-s}\rightarrow \infty$ as $s\rightarrow 1$ which
implies that there are infinitely many primes of the form $a+bx$
when $\gcd (a,b)=1$.

\vspace{3mm}More generally, denote the number of some special kinds
of primes not exceeding $x$ by $P(x)$. If we can prove
$P(x)\rightarrow \infty$ as $x\rightarrow \infty$, then we not only
know the existence of these special kinds of primes, but also know
their infinitude. This is a good method which goes back to Legendre
who firstly conjectured $\pi (x)\approx \frac{x}{\log x -1.08...}$.
Gauss found that a good approximation to $\pi (x)$ is
$li(x)=\int_2^x\frac{dt}{\log t}$. It is easy to prove that $\pi
(x)\geq \frac{\log x}{2\log 2}$ which implies the existence of
infinitely many primes again. In 1851, Tchebychev proved firstly
that for all sufficiently large $x$,  $0.92\frac{x}{\log x}<\pi
(x)<1.10\frac{\ x}{\log x}$. In 1896, Hadamard and de la Vall\'{e}e
Poussin proved independently $\pi (x)\sim \frac{x}{\log x}$ (or
equivalently, $\pi (x)\sim li(x)$ ). This is famous Prime Number
Theorem which implies simply the existence of infinitely many
primes.

 \vspace{3mm}As we mentioned, by studying the behavior
of $P(x)$, one not only can determine the existence and infinitude
of some special kinds of primes, but also know the distribution of
these special kinds of primes, this is a quantitative form and
becomes then a main method for studying the infinitude of some
special kinds of primes. However, it also is the most difficult.
Next section, we would like to focus our attention on the natural
necessary condition and try to give a new sufficient condition of
the infinitude of some special kinds of primes. This leads to a new
way for determining the existence of these primes.

%%%%%%%%%%%%%%%%%%%%%%%%%%%%%%%%%%%%%%%%%%%%%%%%%%%%%%%%%%%%%%%%
\section{A sufficient condition that
a multivariable number-theoretic function represents primes for
infinitely many integral points }
%%%%%%%%%%%%%%%%%%%%%%%%%%%%%%%%%%%%%%%%%%%%%%%%%%%%%%%%%%%%%%%%
In this section, we begin with Euclid's proof of the infinitude of
primes. Euclid's beautiful proof by contradiction goes as follows:
Suppose that there are only finitely many primes, say $k$ of them,
which denoted by $2=p_1<...<p_k$. Note that $p_1...p_k+1 >1$ and
hence it must have a prime factor which differs from $p_1,...,p_k$
and this leads to a contradiction.

\vspace{3mm}Euclid's proof is essentially to construct a number $x$
such that $x$ is coprime to the product  $p_1...p_k$. Note that 2
and 3 are prime. So $|Z_{p_1...p_k}^*|>1$, by Euler function
formula. On the other hand, as we know, if $a$ is the smallest
integer such that $a>1$ and $\gcd (a,p_1...p_k)=1$ then $a$ is
prime. Therefore, there are infinitely many primes since
$|Z_{p_1...p_k}^*|>1$ implies that there is such an integer $a$ in
$Z_{p_1...p_k}^*$. This gives a proof for the infinitude of primes.
Although the proof perhaps is not new, it is enlightened us. This
proof need not construct a new number  $x$  such that  $x$ is
coprime to the product $p_1...p_k$ but prove directly that there is
a number $x>1$ such that $x$ is coprime to the product $p_1...p_k$.
Hence  $x$ has a new prime factor and it leads to a contradiction.
By the existence of such a $x$, there must be the least positive
integer $x>1$ which is coprime to the product  $p_1...p_k$. Of
course, it is prime.

\vspace{3mm}The question of existence of infinitely many primes can
be regard as the question of existence of infinitely many prime
values of the polynomial $f(x)=x$. For any positive integer $m>1$,
$f(S_f(m))=S_f(m)$ always is prime when $f(x)=x$, where $S_f(m)$ is
the least positive integer $n$ such that $\gcd (f(n),m)=1$ and
$f(n)>1$. More generally, let $f(x)$ be a generic number-theoretic
function, unfortunately, for any positive integer $m>1$, $f(S_f(m))$
is not always prime. For example, let $f(x)=2^x-1$ and $m=82677$,
$S_f(m)=11$ and $f(S_f(m))=2^{11}-1=23\times 89$ is not prime. Thus
a key fact which states that if $a>1$ is the smallest integer such
that $\gcd (a,m)=1$ and then $a$ is prime is not true in the generic
case. Why is it a key fact that if $a>1$ is the smallest integer
such that $\gcd (a,m)=1$ and then $a$ is prime. As we know, if the
number of primes is finite, then the proposition which states if
$a>1$ is the smallest integer such that $\gcd (a,m)=1$ and then $a$
is prime is false. Therefore, we want to use this fact.
Unfortunately, in the generic case, it is not always true. How to
treat with it?

\vspace{3mm}Let's look back Euclid's proof again. He considered the
product of primes $p_1...p_k$. Similarly, we may consider $p_k!$. In
fact, $p_k!+1$ and $p_1...p_k$ are coprime, which implies the
infinitude of primes again. Directly or more expediently, we
consider the factorial $n!$ instead of the finite product
$p_1...p_k$ of primes. Clearly, so long as $n>p_k$, then it will
lead to a contradiction still. Particularly, let $a\in Z_{n!}^*$ be
the smallest integer such that $a>1$ and $\gcd (a,n!)=1$, then $a$
is prime. This is a key fact. We hope naturally this key fact still
is true in the generic case that a number-theoretic function $f(x)$
or  $f(x_1,...,x_k)$ represents infinitely many primes.  In the
following conjecture 3, we try to give a primary consideration.

\vspace{3mm}Another reason that we would like to consider the
factorial is because the factorial can be viewed as a special case
of the $\Gamma$ function which is closely related to the
distinguished Riemann Hypothesis.

\vspace{3mm}Below is the third reason that we would like to consider
the factorial:

\vspace{3mm} We notice that if a number-theoretic function $f(x)$
represents primes for infinitely many natural numbers $x$, then for
any positive integer $n$, there is a natural numbers $x$ such that
the least prime divisor of $f(x)$ is greater than  $n$. Therefore,
there must be a least natural numbers $k$ such that the least prime
divisor of $f(k)$ is greater than  $n$. Namely, $f(k)$ ($>1$) is
coprime to  $n!$. We also know that there must be a least natural
numbers $r$ such that  $f(r)$ ($>1$) is coprime to  $n!$. Of source,
$r=k$. Very naturally, one might believe that $f(k)=f(r)$ is prime.

\vspace{3mm}The following Proposition 3 further gives some
witnesses.

\vspace{3mm}\noindent {\bf Proposition 3:~~}%
Let $f(x)$ be a generic number-theoretic function. If there is a
constant $c$ such that for every positive integer $m>c$, there is a
natural number $y$ such that $f(y)>1$ is in $Z_m^* $, and if the
least number $f(x)$ which exceeds 1 in $Z_m^* $ is not prime, then
$f(x)$ represents primes at most for finitely many natural numbers
$x$.

\vspace{3mm}\noindent {\bf Proof of Proposition 3:~~}%
If $f(x)$ represents primes for infinitely many natural numbers $x$,
then there is a natural number $y$ such that $f(y)>c$ is prime.
Without loss of generality, assume that $f(y)$ is the least prime
which exceeds $c$. If $f(y)>2$, then $2((f(y)-1)!)>c$. But $f(y)$ is
prime and also is the least natural number which exceeds 1 in
$Z_{2((f(y)-1)!)}^* $. By our assumption, this least natural number
 should be a composite number. This is a contradiction. Therefore,
 $f(y)=2$ and $2>c$. In this case, note
 that $3>c$
and we have $f(y)=2\in Z_3^*$. But, $2$ is prime and also is the
least number which exceeds 1 in $Z_3^*$. This is a contradiction
again. Therefore, Proposition 3 holds.

\vspace{3mm}By the proof of this proposition, we see also that if
$f(x)$ represents primes for infinitely many numbers $x$, then there
are infinitely many numbers $m$ such that the least number $f(y)$
which exceeds 1 in $Z_m^*$ is prime. One could generalize it to the
generic case. We also believe naively that if $f(x)$ represents
primes for infinitely many numbers $x$, then there is a positive
integer $c$ such that for each $m>c$,  if $f(r)$ ($>1$) is the least
natural numbers  of the form $f(r)$ such that  $f(r)$ is coprime to
$m!$, then $f(r)$ is prime.

\vspace{3mm} Due to the fact the $f(x)=ax+b$ with $\gcd (a,b)=1$
represents primes for infinitely many natural numbers $x$, we now
prove that there is a positive integer $c$ such that for each $m>c$,
and if $f(r)$ is the least prime of the form $f(x)=ax+b$ such that
$\gcd (f(r), m!)=1$, then $f(r)< m!$. This is easy to prove. Denote
the $i^{th}$ prime of the form $f(x)=ax+b$ by $P_{f,i}$. In [103],
we have proved that there is a constant $C$ depending on $a$ and $b$
such that when $n>C$, $\prod_{i=1}^{i=n}P_{f,i}>P_{f,n+1}$. Let
$k\geq n$ and
$\prod_{i=1}^{i=k}P_{f,i}<m!<\prod_{i=1}^{i=k+1}P_{f,i}$. Clearly,
$\gcd (P_{f,k+1},m!)=1$. So, $P_{f,k+1}\geq f(r)$. If $f(r)\geq m!$,
then $f(r)\geq m!>\prod_{i=1}^{i=k}P_{f,i}>P_{f,k+1}$. It is a
contradiction.  Thus $f(r)< m!$. By the results in [103],  we still
have similar results for the cases $f(x,y,z,w)=x^2+y^2+z^2+w^2$,
$f(x,y)=x^2+y^2+1$, $f(x,y)=x^2+y^4$, $f(x,y)=x^3+2y^3$, and so on.

\vspace{3mm}Based on the discussion above, now, we give a sufficient
condition of Conjecture 2 as follows.

\vspace{3mm}\noindent {\bf Conjecture 3:~~}%
Let $f(x_1,...,x_k)$ be a multivariable polynomial with  integral
coefficients (or a multivariable number-theoretic function), if
there is a positive integer $c$ such that for every positive integer
$m\geq c$, there exists an integral point $(y_1,...,y_k)$ such that
$f(y_1,...,y_k)>1$ is in $Z_m^* $, and there exists an integral
point $(z_1,...,z_k)$ such that $f(z_1,...,z_k)\geq c$ is primes,
moreover, for any integer $l$ with $l\geq r$, if $f(x_1,...,x_k)$ is
the least positive inetger such that $f(x_1,...,x_k)>1$ is in
$Z_{l!}^* $, then $f(x_1,...,x_k)$ represents  primes, where $r$  is
the least positive integer such that provided $n\geq r!$, then there
exists an integral point $(y_1,...,y_k)$ such that
$f(y_1,...,y_k)>1$ is in $Z_n^* $.

\vspace{3mm}As for more generic case of several multivariable
number-theoretic functions, it is very complicated. For instance,
let $f_1(x)=x,f_2(x)=x+180$. It is easy to prove that for each
$n>5$, there is a least natural numbers $x$ such that
$f_1(x)=x>1,f_2(x)=x+180>1$ and $f_1(x)\times f_2(x)\in Z_{n!}^*$.
But when $n=6$, we have $f_1(x)=x=7$ and $f_2(x)=187$ is not prime.
We left this question to the readers. However, when
$f_1(x_1,...,x_k),...,f_s(x_1,...,x_k)$ are multivariable
polynomials with  integral coefficients, we fix  Conjecture 3 and
further generalize it as follows:

\vspace{3mm} For $s>1$, let $f_1(x_1,...,x_k),...,f_s(x_1,...,x_k)$
be multivariable polynomials with  integral coefficients, assume
that there is a positive integer $c$ such that for every positive
integer $m\geq c$, there exists an integral point $(y_1,...,y_k)$
such that $f_1(y_1,...,y_k)>1,...,f_s(y_1,...,y_k)>1$ are all in
$Z_m^* $, and there exists an integral point $(z_1,...,z_k)$ such
that $f_1(z_1,...,z_k)\geq c,...,f_s(z_1,...,z_k)\geq c$ are all
primes. Then for any integer $l$ with $l\geq r$, there exists an
integral point $(x_1,...,x_k)$ such that
$f_1(x_1,...,x_k),...,f_s(x_1,...,x_k)$ are all in $Z_{l!}^* $ and
$f_1(x_1,...,x_k),...,f_s(x_1,...,x_k)$  represent simultaneously
primes, where $r$  is the least positive integer such that provided
$n\geq r!$, then there exists an integral point $(w_1,...,w_k)$ such
that $f_1(w_1,...,w_k)>1,...,f_s(w_1,...,w_k)>1$  are all in $Z_n^*
$.

%%%%%%%%%%%%%%%%%%%%%%%%%%%%%%%%%%%%%%%%%%%%%%%%%%%%%%%%%%%%%%%%
\section{A sufficient and necessary
condition that several multivariable number-theoretic functions
represent simultaneously primes for infinitely many integral points}
%%%%%%%%%%%%%%%%%%%%%%%%%%%%%%%%%%%%%%%%%%%%%%%%%%%%%%%%%%%%%%%%
After finishing Section 8, the author felt intensively that there
must be a sufficient and necessary condition that several
multivariable number-theoretic functions represent simultaneously
primes for infinitely many integral points. Therefore, this section
was added very recently. The author obtrusively suggested a
generalization of Conjecture 2 as follows.

\vspace{3mm}\noindent {\bf  A sufficient and necessary
condition:~~}%
Let $f_1(x_1,...,x_k),...,f_s(x_1,...,x_k)$ be multivariable
number-theoretic functions, assume that there is a positive integer
$c$ such that for every positive integer $m\geq c$, there exists an
integral point $(y_1,...,y_k)$ such that
$f_1(y_1,...,y_k)>1,...,f_s(y_1,...,y_k)>1$ are all in $Z_{m!}^* $,
and there exists an integral point $(z_1,...,z_k)$ such that
$f_1(z_1,...,z_k)\geq c!,...,f_s(z_1,...,z_k)\geq c!$ are all
primes. Then the sufficient and necessary condition that
$f_1(x_1,...,x_k)$, ... , $f_s(x_1,...,x_k)$  represent
simultaneously primes for infinitely many integral points is  that
for every $m$, there is an integral point $(m_1,...,m_k)$ such that
 $f_1(m_1,...,m_k)$, ... , $f_s(m_1,...,m_k)$ in $Z_{m!}^*$ are all primes.

\vspace{3mm} Due to the fact that we consider the factorial which
can be viewed as a special case of the $\Gamma$ function, we should
assume that  $f_1,...,f_s$ are continuous functions on $R^k$, where
$R$ is the set of all real numbers.

\vspace{3mm} This sufficient and necessary condition implies that
there is a positive integer $c$ such that for every positive integer
$m\geq c$, there exists an integral point $(y_1,...,y_k)$ such that
$f_1(y_1,...,y_k)>1,...,f_s(y_1,...,y_k)>1$ are all in $Z_{m!}^* $.
It looks slightly weaker than Conjecture 1. Thus, one could ask:
does it imply that there are only finitely many Fermat primes?
Sensuously, it seems that $f(x)=2^{2^x}+1$ satisfies this weakened
necessary condition. Namely,  there is a positive integer $c$
(perhaps $c=3$) such that for every positive integer $m\geq c$,
there exists an integer $x>0$ such that $2^{2^x}+1$ is in $Z_{m!}^*
$. But, if $2^{2^a}+1$ is the least Fermat number in $Z_{m!}^* $,
then, $2^{2^a}+1$ might is not always prime. For instance, by
Stirling's formula $n!\approx \sqrt{2\pi n}(\frac{n}{e})^n$ and the
factorization of Fermat numbers, we have $2^{2^6}+1\in
Z_{(2^{2^4}+1)!}^* $, $2^{2^7}+1\in Z_{(2^{2^5}+1)!}^* $ and so on,
but $2^{2^6}+1$ and $2^{2^7}+1$ are the least but not prime
respectively. Does it imply that  there are only finitely many
Fermat primes again?

\vspace{3mm}Anyway, are these conjectures proposed consistent with
each other? Are they reasonable or reliable? The author are waiting
for advice of readers. With the development of mathematics, the
correct answers will come---we must know, we will know, Hilbert
said.

%%%%%%%%%%%%%%%%%%%%%%%%%%%%%%%%%%%%%%%%%%%%%%%%%%%%%%%%%%%%%%%%
\section{Conclusion}
%%%%%%%%%%%%%%%%%%%%%%%%%%%%%%%%%%%%%%%%%%%%%%%%%%%%%%%%%%%%%%%%
I learn from Euclid all the time. This paper is a part of my paper
\emph{Euclid's algorithm and the infinitude of some special kinds of
primes}, in which, his two great number-theoretical
achievements---Euclid's algorithm and his proof for the infinitude
of primes, have been studied. Of course, these two significant
results are not independent, and, Chinese Remainder Theorem is a
bridge between them because Euclid's algorithm implies Chinese
Remainder Theorem which also implies Euclid's second theorem.

\vspace{3mm}On Euclid's algorithm, we have done the following work
in the paper \emph{Euclid's algorithm and the infinitude of some
special kinds of primes}: (1)Euclid's Number-Theoretical Work [77];
(2)Euclid's Algorithm, Guass' Elimination and Buchberger's Algorithm
[85]; (3)Euclid's Algorithm and W Sequences [86]; (4)Euclid's
Algorithm and three public key cryptosystems---RSA Cryptosystem,
Elliptic Curve Cryptosystems and Multivariate Public Key
Cryptosystems; (5) Euclid's Algorithm, LLL Algorithm and the number
field sieve. On Euclid's proof for the infinitude of primes, it
leads to this paper. Knuth [87] called Euclid's Algorithm the
granddaddy of all algorithms. In the author's eyes, Euclid's proof
for the infinitude of primes also is the granddaddy of some proofs
for the infinitude of some kind special kinds of primes.

\vspace{3mm}In this paper, we try to establish a generic model for
the problem of infinitude of some special kinds of primes. More
precisely, we try to establish a generic model for the problem that
several multivariable number-theoretic functions represent
simultaneously primes for infinitely many integral points. We
analyzed some equivalent necessary conditions that irreducible
univariable polynomials with integral coefficients represent
infinitely many primes, found new necessary conditions which perhaps
imply that there are only finitely many Fermat primes, generalized
Euler's function, the prime-counting function and
Schinzel-Sierpinski's Conjecture and so on, obtained an analogy of
the Chinese Remainder Theorem. Finally, a sufficient and necessary
condition that several multivariable number-theoretic functions
represent simultaneously primes for infinitely many integral points
was proposed. Nevertheless, this is only a beginning and it miles to
go. The author would like to cite the comment in Schinzel and
Sierpinski's paper "\emph{we do not know what will be the fate of
our hypothesis, however, we think that, even if they are refuted,
this will not be without profit for Number Theory.}" to close this
paper. Please let me know any questions, reviews and criticisms
at\emph{shaohuazhang@mail.sdu.edu.cn}. Thank you very much.

%%%%%%%%%%%%%%%%%%%%%%%%%%%%%%%%%%%%%%%%%%%%%%%%%%%%%%%%%%%%%%%%
\section{Acknowledgements}
%%%%%%%%%%%%%%%%%%%%%%%%%%%%%%%%%%%%%%%%%%%%%%%%%%%%%%%%%%%%%%%%
God created the integers. Thank God for the great blessings he has
given me to study integers.

\vspace{3mm}God could not be everywhere and therefore he made
mothers. Thank my mother for arousing my interest in natural
numbers. She used to tell me: "seven is sacred; two is female; three
is male; and odd numbers are not perfect; and so on". Even now, I
believe that there are no odd perfect numbers. Unfortunately, she
died soon after I began to read the number-theoretical book for the
first time. From then on, I often bring to mind her smile and her
encouragement when I come into contact with integers. Today, I
dedicate this paper to her. I regret that I have not proved an
excellent number-theoretical theorem for her. What I can do now is
to try my best to become a disseminator of some excellent
number-theoretical theorems. In [Appendix], I give a list of 100
theorems in Number Theory.

%%%%%%%%%%%%%%%%%%%%%%%%%%%%%%%%%%%%%%%%%%%%%%%%%%%%%%%%%%%%%%%%
\section{Appendix}
%%%%%%%%%%%%%%%%%%%%%%%%%%%%%%%%%%%%%%%%%%%%%%%%%%%%%%%%%%%%%%%%
{\bf $$\mbox {100 theorems in Number Theory}$$ ~~} "Number Theory is
the queen of mathematics"---Guass. For a long time, I want to edit a
number-theoretical e-book which includes many excellent theorems and
is worthwhile to take my lifetime to learn. Recently, this work is
almost completed. I named it "100 theorems in Number Theory". If
Number Theory is the queen of mathematics, then these theorems are
her pearls. Of course, there are much more than 100 theorems in
Number Theory. To follow a principle that Mathematics is essentially
simple, and based on the individual opinion and taste, I only pick
some theorems which are my favorites and look simple in spite of
some proofs are extremely intricate. Namely, the description of
these theorems is very easy, although their proofs maybe be
extremely difficult. For example, the meaning of Fermat's last
theorem or Green-Tao theorem  is very clear, but its proof is
difficult to understand. Unfortunately, the proofs of about 51
theorems in this e-book, while intensely enjoyable, do require hard
study to grasp. Maybe, this is a basic reading for Number Theory. If
one is afraid of meeting the difficulty, then he always meets the
difficulty. Therefore, anyone who loves Number Theory should learn
the proofs of approximately 87 theorems. I will try my best to
travel for this dream. A Bachelor of Number Theory had better
understand the proofs of more or less 60 theorems. A Master of
Number Theory had better understand the proofs of more or less 70
theorems. And a Doctor of Number Theory had better understand the
proofs of more or less 80 theorems. Below is the list of theorems.\\

\vspace{3mm}\noindent {\bf 1 The first theorem about Theory of
Divisibility
(also called division algorithm):~~}%
Let $a$ and  $b$ be integers with  $b>0$ . There exist unique
integers $q$  and  $r$ such that $a=bq+r$ and  $0\leq r<b$ .

\vspace{3mm}\noindent {\bf Remark:~~}%
This theorem is the basis of Theory of Divisibility. Many
number-theoretical texts begin with it. However, Euclid did not do
like this. Euclid began his number-theoretical work by introducing
his algorithm which states essentially that for two distinct
positive integers, replace continually the larger number by the
difference of them until both are equal, then the answer is their
greatest common divisor. In [77], we showed that Euclid's algorithm
is equivalent with Division algorithm.

\vspace{3mm}\noindent {\bf  2 Euclid's first theorem[67]:~~}%
 If $p$ is prime, and $p|ab$, then $p|a$ or $p|b$.

\vspace{3mm}\noindent {\bf 3 Euclid's second theorem[67]:~~}%
 The number of primes is infinite.

\vspace{3mm}\noindent  {\bf 4 The  fundamental theorem of arithmetic:~~}%
Every positive integer can be written uniquely (up to order) as the
product of prime numbers.

\vspace{3mm}\noindent  {\bf 5 The linear congruence theorem:~~}%
If $a$ and $b$ are any integers and $n$ is a positive integer, then
the congruence $ax\equiv b (\mod n)$ has a solution for $x$ if and
only if $b$ is divisible by the greatest common divisor $(a,n)$ of
$a$ and $n$. Particularly, when $(a,n)=1$, the congruence $ax\equiv
b (\mod n)$ has a unique solution modulo $n$.

\vspace{3mm}\noindent  {\bf 6 Chinese Remainder Theorem:~~}%
Given a system of simultaneous linear congruences  $x\equiv a_i
(\mod n_i)$ for $i=1,2,...,k$ and for which $n_i$ are pairwise
relatively prime positive integers, where $a_i$ are integers, then
this linear system has a unique solution modulo
$n=\prod_{i=1}^{i=k}n_i$. Particularly, for $i=1,2,...,k$, if $\gcd
(a_i,n_i)=1$, then, this linear system has a unique solution $x$ in
$Z_n^*$.

\vspace{3mm}\noindent {\bf 7 Fermat's little theorem:~~}%
Let $p$ be a prime, if the integer $a$ is not divisible by $p$, then
$a^{p-1}\equiv 1 (\mod p)$. Moreover, $a^p\equiv a(\mod )p$ for
every integer $a$.

\vspace{3mm}\noindent {\bf 8 Euler's theorem:~~}%
Let $m$ be a positive integer, and let $a$ be an integer relatively
prime to $m$. Then $a^{\varphi (m)}\equiv 1 (\mod m)$, where
$\varphi (m)$ is defined to be the number of positive integers less
than or equal to $m$ that are coprime to $m$.

\vspace{3mm}\noindent {\bf 9 Carmichael's theorem:~~}%
If $(a,m)=1$, then $a^\lambda \equiv 1 (\mod m)$, where $\lambda$ is
the smallest integer such that $k^\lambda \equiv 1 (\mod m)$ for all
$k$ relatively prime to $m$.

\vspace{3mm}\noindent {\bf 10 Wilson's theorem:~~}%
If $p$ is prime, then $(p-1)! \equiv -1 (\mod p)$.

\vspace{3mm}\noindent {\bf 11 A theorem of Wolstenholme:~~}%
If $p>3$ is prime, then the numerator of harmonic number
$1+\frac{1}{2}+...+\frac{1}{p-1}$ is divisible by $p^2$.

\vspace{3mm}\noindent {\bf 12 Euclid-Euler theorem:~~}%
$n$  is even perfect number if and only if $n=2^{p-1}(2^p-1)$, where
$p$ is prime such that $2^p-1$ is also prime.

\vspace{3mm}\noindent {\bf 13 Lam\'{e}'s theorem:~~}%
Finding the greatest common divisor of integers $a$ and $b$ with
$a>b$, Euclid's algorithm runs in no more than $5k$ steps, where $k$
is the number of (decimal) digits of $b$.

\vspace{3mm}\noindent {\bf 14 Midy's theorem:~~}%
If the period of a repeating decimal for $\frac{a}{p}$, where $p$ is
prime and $\frac{a}{p}$ is a reduced fraction, has an even number of
digits, then dividing the repeating portion into halves and adding
gives a string of 9s. For example,
$\frac{1}{7}=0.\overline{142857}$, $142+857=999$.

\vspace{3mm}\noindent {\bf 15 Bauer's theorem:~~}%
Let $m>2$ be an integer and let $f(x)$ be an integral polynomial
that has at least one real root. Then $f(x)$ has infinitely many
prime divisors that are not congruent to $1(\mod m)$.

\vspace{3mm}\noindent {\bf 16 Euler's quadratic residue theorem:~~}%
Let $p$ be an odd prime. For every integer $a$,
 $(\frac{a}{p})\equiv a^{\frac{p-1}{2}}(\mod p)$, where $(\frac{a}{p})$ is the Legendre symbol.

\vspace{3mm}\noindent {\bf 17 The golden theorem (also called the
law of quadratic
reciprocity):~~}%
If $p$ and $q$ is distinct odd primes, then
$(\frac{q}{p})(\frac{p}{q})=(-1)^{(\frac{p-1}{2})(\frac{q-1}{2})}$.

\vspace{3mm}\noindent {\bf 18 Fermat's theorem on sums of two squares:~~}%
An odd prime $p$ is expressible as $p=x^2+y^2$ with $x$ and $y$ are
integers, if and only if $p\equiv 1(\mod 4)$.

\vspace{3mm}\noindent {\bf 19 Lagrange's four-square theorem:~~}%
Every positive integer can be expressed as the sum of four squares
of integers.

\vspace{3mm}\noindent {\bf 20 Fermat polygonal number theorem:~~}%
Every positive integer is a sum of at most $n$  n-polygonal numbers,
where $n>2$  is a positive integer.

\vspace{3mm}\noindent {\bf 21 A theorem of Carmichael on the n-th Fibonacci number:~~}%
Every Fibonacci number $f_n$ with $n\neq 1, 2, 6, 12$, has at least
one characteristic factor which is not a factor of any earlier
Fibonacci number.

\vspace{3mm}\noindent {\bf 22 Lagrange's continued fraction theorem:~~}%
A number is a quadratic surd if and only if its continued fraction
expansion is eventually periodic.

\vspace{3mm}\noindent {\bf 23 Dirichlet's approximation theorem:~~}%
Given any real number $\theta$ and any positive integer $n$, there
exist integers $x$ and $y$ with $0<x\leq n$ such that $|x\theta
-y|<\frac{1}{n}$.

\vspace{3mm}\noindent {\bf 24 Liouville's theorem:~~}%
For any algebraic number $x$ with degree $n>1$, there exists $c>0$
such that  $|x-\frac{p}{q}|>\frac{c}{q^n}$ for all rationals
$\frac{p}{q}(q>0)$.

\vspace{3mm}\noindent {\bf 25 Van der Waerden's theorem:~~}%
If the set of all positive integers is written as the union of sets
of the finite number, then there exists at least a set which
contains arbitrarily long arithmetic progressions.

\vspace{3mm}\noindent {\bf 26 Minkowski's theorem:~~}%
Any convex set in $R^n$ which is symmetric with respect to the
origin and with volume greater than $2^n$ contains a non-zero
lattice point.

\vspace{3mm}\noindent {\bf 27 Bertrand-Chebyshev theorem:~~}%
There exists a prime in interval $(n,2n)$ when $n>1$.

\vspace{3mm}\noindent {\bf 28 Mills' theorem:~~}%
There exists a real constant $\theta$ such that $[\theta^{3^n}]$ is
prime for all $n\geq 1$.

\vspace{3mm}\noindent {\bf 29 Rosser's theorem:~~}%
Let $p_n$ be the n-th prime number, then for $n>1$, $p_n>n\ln n$.

\vspace{3mm}\noindent {\bf 30 Beatty's theorem:~~}%
Let $x$ and $y$ be positive irrational numbers satisfying
$\frac{1}{x}+\frac{1}{y}=1$. Then each positive integer belongs to
exactly one of the two sequences $\{[nx]\}$ and $\{[ny]\}$.

\vspace{3mm}\noindent {\bf 31 Hurwitz's Irrational Number Theorem:~~}%
For any irrational number $\theta$ , there are infinitely many
rational numbers $\frac{x}{y}$ such that $|\theta
-\frac{x}{y}|<\frac{1}{\sqrt{5}y^2} $.

\vspace{3mm}\noindent {\bf 32 Blichfeldt's Theorem:~~}%
A bounded set of points $C$ with area $A$, can be translated to a
position $C'$ so as to cover a number of lattice points greater than
$A$.

\vspace{3mm}\noindent {\bf 33 Ramanujan-Skolem's theorem:~~}%
The equation $x^2+7=2^n$ has solutions in natural numbers $n$ and
$x$ just when $n=3,4,5,7,15$.

\vspace{3mm}\noindent {\bf 34 Thue's Theorem:~~}%
If $f$ is a bivariate form with rational coefficients which is
irreducible over the rational numbers and has degree $\geq 3$, and
$r$ is a rational number other than 0, then the equation $f(x,y)=r$
has only finitely many solutions in integers $x$ and $y$.

\vspace{3mm}\noindent {\bf 35 Rotkiewicz Theorem:~~}%
If $n\geq 19$, there exists a poulet number between $n$ and $n^2$.

\vspace{3mm}\noindent {\bf Remark:~~}%
A poulet number $m$ is a Fermat pseudoprime to base 2, namely, a
composite number $m$ satisfying $2^{m-1}\equiv 1 (\mod m)$.

\vspace{3mm}\noindent {\bf 36 Schnirelmann's Theorem:~~}%
There exists a positive integer C such that every sufficiently large
integer is the sum of at most C primes.

\vspace{3mm}\noindent {\bf 37 Sierpinski's Composite Number Theorem:~~}%
There exist infinitely many positive odd numbers $k$ such that
$k2^n+1$ is composite for every integer $n>0$.

\vspace{3mm}\noindent {\bf 38 Sierpinski's Prime Sequence Theorem:~~}%
For any $m$, there exists a number $k$ such that the sequence
$\{n^2+k\}$ contains at least $m$ primes.

\vspace{3mm}\noindent {\bf 39 A theorem of H.Gupta and S.P. Khare:~~}%
For $2<k<1794$ , $\binom{k^2}{k}$ is greater than the product of the
first $k$ primes, while for $k\geq 1794$, $\binom{k^2}{k}$ is less
than the product of the first $k$ primes.

\vspace{3mm}\noindent {\bf 40 A theorem of Sylvester and Schur:~~}%
The product of $k$ consecutive positive integers each exceeding $k$
is divisible by a prime greater than $k$.

\vspace{3mm}\noindent {\bf 41 Theorem of Pillai and Szekeres:~~}%
For any positive integer $n\geq 17$, there exists a sequence of $n$
consecutive positive integers such that no one of this sequence is
relatively prime with all of the others.

\vspace{3mm}\noindent {\bf 42 Erd\"{o}s-Anning Theorem:~~}%
An infinite number of points in the plane can have mutual integer
distances only if all the points lie on a straight line.

\vspace{3mm}\noindent {\bf 43 A theorem of Erd\"{o}s and Selfridge:~~}%
The product of consecutive integers is never a power.

\vspace{3mm}\noindent {\bf 44 The theorem of Waring-Hilbert:~~}%
Every positive integer $n$ is the sum of at most s k-th powers of
natural numbers, where  $s=s(k)$ is independent of $n$.

\vspace{3mm}\noindent {\bf 45 Mason's theorem:~~}%
Let $f,g,h$ be three polynomials with no common factors such that
$f+g=h$. Then the number of distinct roots of the three polynomials
is either one or greater than their largest degree.

\vspace{3mm}\noindent {\bf 46 Dirichlet's theorem:~~}%
For any two positive coprime integers $a$ and $b$, there are
infinitely many primes of the form $a+bn$.

\vspace{3mm}\noindent {\bf 47 Chebyshev's theorem:~~}%
For real number $x$, denote the number of primes less than or equal
to $x$ by $\pi (x)$. Then there exist positive constants  $A$ and
$B$ such that $Ax\leq \pi (x) \ln x \leq Bx$.

\vspace{3mm}\noindent {\bf 48 Prime number theorem:~~}%
$\pi (x)\sim \frac{x}{\ln x}$.

\vspace{3mm}\noindent {\bf 49 A theorem of G. Robin:~~}%
Denote the number of distinct prime factors of $x$ by $\omega (x)$,
then for every integer $n\geq 26$, $\omega (n)<\frac{\log
n}{\log\log n -1.1714}$, with equality when $n$ is the product of
the first 189 primes.

\vspace{3mm}\noindent {\bf 50 A theorem of M. Agrawal -N. Kayal and N. Saxena:~~}%
There is an algorithm determines whether a number is prime or
composite within polynomial time.

\vspace{3mm}\noindent {\bf 51 Brun's theorem:~~}%
The sum of the reciprocals of the twin primes is convergent with a
finite value.

\vspace{3mm}\noindent {\bf 52 Ap\'{e}ry's theorem:~~}%
The number $\zeta (3)=1+\frac{1}{2^3}+\frac{1}{3^3}+...+...$ is
irrational.

\vspace{3mm}\noindent {\bf 53 A theorem of Motohashi Y.:~~}%
There are infinitely many primes of the form $x^2+y^2+1$.

\vspace{3mm}\noindent {\bf 54 A theorem of
Fouvry, Etienne and Iwaniec H.:~~}%
There are infinitely many primes of type $x^2+y^2$, where $x$ is a
prime number.

\vspace{3mm}\noindent {\bf 55 A theorem of Friedlander, John and Iwaniec H.:~~}%
There are infinitely many primes of type  $x^2+y^4$.

\vspace{3mm}\noindent {\bf 56 A theorem of Heath-Brown:~~}%
There are infinitely many primes of form  $x^3+2y^3$.

\vspace{3mm}\noindent {\bf 57 Linnik's theorem:~~}%
Denote $p(l,k)$ the least prime in the arithmetic progression $l+kn$
with $(l,k)=1$, where $n$ runs through the positive integers, and
let $p(k)$ be the maximum value of $p(l,k)$  for all $l$ satisfying
$(l,k)=1$ and $1\leq l\leq k$. Then there exist positive $C$ and $L$
such that $p(k)<Ck^L$.

\vspace{3mm}\noindent {\bf 58 Heath-Brown's theorem:~~}%
(with the notation above) $p(k)<Ck^{5.5}$.

\vspace{3mm}\noindent {\bf 59 Vinogradov's theorem:~~}%
Every sufficiently large odd number can be written as the sum of
three primes.

\vspace{3mm}\noindent {\bf 60 Chen's theorem:~~}%
Every sufficiently large even number can be written as the sum of
either two primes, or a prime and a semiprime.

\vspace{3mm}\noindent {\bf 61 A theorem of Roth:~~}%
For every value $d$ with $0<d<1$, there is a number $C$ such that
every subset $A$ of $\{1,2,3,...,N\}$ of cardinality $dN$ contains a
length-3 arithmetic progression, provided $N>C$.

\vspace{3mm}\noindent {\bf 62 Szemeredi's theorem:~~}%
Every sequence of integers that has positive upper density contains
arbitrarily long arithmetic progressions.

\vspace{3mm}\noindent {\bf 63 A theorem of J. G. van der Corput:~~}%
The primes contain infinitely many arithmetic progressions of length
3.

\vspace{3mm}\noindent {\bf 64 A theorem of Ben Green:~~}%
Any set containing a positive proportion of the primes contains a
3-term arithmetic progression.

\vspace{3mm}\noindent {\bf 65 A theorem of Balog:~~}%
For any $m>1$, there are $m$ distinct primes $p_1,...,p_m$  such
that all of the averages $\frac{p_i+p_j}{2}$ are primes.

\vspace{3mm}\noindent {\bf 66 Green-Tao theorem:~~}%
The primes contain arbitrarily long arithmetic progressions.

\vspace{3mm}\noindent {\bf 67 A theorem of W. R. Alford -A. Granville and C. Pomerance:~~}%
There are infinitely many Carmichael numbers.

\vspace{3mm}\noindent {\bf 68 Tijdeman's theorem:~~}%
There are at most a finite number of consecutive powers.

\vspace{3mm}\noindent {\bf 69 Mihailescu's theorem:~~}%
8 and 9 are the only consecutive powers.

\vspace{3mm}\noindent {\bf 70 A theorem of Pythagoras' school:~~}%
$\sqrt{2}$ is irrational.

\vspace{3mm}\noindent {\bf 71 A theorem of Euler on the irrationality:~~}%
The base of the natural logarithm $e$ is irrational.

\vspace{3mm}\noindent {\bf 72 A theorem of Lambert:~~}%
The ratio $\pi$ of a circle's circumference to its diameter is
irrational.

\vspace{3mm}\noindent {\bf 73 Hermite-Lindemann theorem:~~}%
$e$  and $\pi$ are all transcendental numbers.

\vspace{3mm}\noindent {\bf 74 Gelfond-Schneider theorem:~~}%
If $\alpha$ and  $\beta$ are algebraic numbers (with $\alpha \neq 0,
1$), and if  $\beta$ is not a rational number, then any value
 of $\alpha^\beta$ is a transcendental number.

\vspace{3mm}\noindent {\bf 75 Six Exponentials Theorem:~~}%
Let $(x_1,x_2)$ and $(y_1,y_2,y_3)$ be two sets of complex numbers
linearly independent over the rational number field. Then at least
one of
$e^{x_1y_1},e^{x_1y_2},e^{x_1y_3},e^{x_2y_1},e^{x_2y_2},e^{x_2y_3}$
is transcendental.

\vspace{3mm}\noindent {\bf 76 Baker-Stark theorem:~~}%
The only imaginary quadratic fields $Q(\sqrt{-d})$  with class
number 1, where $d$ is a square-free positive integer, are given by
$d =1, 2, 3, 7, 11, 19, 43, 67, 163$.

\vspace{3mm}\noindent {\bf 77 Thue-Siegel-Roth theorem:~~}%
For any given algebraic number $\theta$, and for given $\varepsilon
>0$, the inequality $|\theta-\frac{p}{q}|<\frac{1}{q^{2+\varepsilon}}$ can have only finitely many solutions in coprime
integers $p$ and $q$.

\vspace{3mm}\noindent {\bf 78 A theorem of Yu. V. Nesterenko:~~}%
$e$, $\pi$ and $\Gamma (\frac{1}{4})$ are algebraically independent.

\vspace{3mm}\noindent {\bf 79 Matiyasevich's theorem:~~}%
Every recursively enumerable set is Diophantine.

\vspace{3mm}\noindent {\bf 80 Fifteen Theorem:~~}%
If an integral quadratic form with integral matrix represents all
positive integers up to 15, then it represents all positive
integers.

\vspace{3mm}\noindent {\bf 81 Erd\"{o}s -Kac theorem:~~}%
If $\omega (n)$ is the number of distinct prime factors of $n$, then
for any fixed $a<b$, $\lim_{N\rightarrow\infty}\frac{1}{N}|\{n\leq
N: a\leq \frac{\omega (n)-\log\log N}{\sqrt{\log\log N}}\leq
b\}|=\int_a^b \varphi (x) dx$, where $\varphi
(x)=\frac{1}{\sqrt{2\pi}}e^{-\frac{x^2}{2}}$ is the probability
density function of the standard normal distribution.

\vspace{3mm}\noindent {\bf 82 A theorem of Pomerance and Selfridge:}
For any given integer $n$ and $m$ with $n>0$, there exists a 1-1
correspondence $f:\{1,...,n\}\rightarrow\{m+1,...,m+n\}$ such that
$\gcd (i,f(i))=1$ for $1\leq i\leq n$.

\vspace{3mm}\noindent {\bf 83 Dirichlet's unit theorem:~~}%
The rank of the group of units in the ring of algebraic integers of
a number field $F$ equals to $r_1+r_2-1$, where $r_1$ is the number
of real embeddings and $r_2$ the number of conjugate pairs of
complex embeddings of $F$.

\vspace{3mm}\noindent {\bf 84 The Fundamental Theorem of Ideal Theory:~~}%
In the domain of all algebraic integers in an algebraic number
field, every nonzero ideal can be represented uniquely (except for
order) as a product of powers of distinct prime ideals.

\vspace{3mm}\noindent {\bf 85 Kronecker-weber theorem:~~}%
Every abelian field is a subfield of a cyclotomic field. Namely, any
Galois extension of the field $Q$ of rational numbers whose Galois
group is Abelian must be a subextension of a field obtained from $Q$
by adjoining root of unity.

\vspace{3mm}\noindent {\bf 86 Kummer's theorem:~~}%
If $p$ is a regular prime which does not divide the class number of
the cyclotomic field $Q(\zeta_p)$, then the equation $x^p+y^p=z^p$
is unsolvable in nature number.

\vspace{3mm}\noindent {\bf 87 Hilbert's basis theorem:~~}%
Every ideal in the ring $F[x_1,...,x_n]$ is finitely generated,
where $F[x_1,...,x_n]$ is a polynomial ring in $n$ variables over a
field $F$.

\vspace{3mm}\noindent {\bf 88 Hilbert's zero theorem:~~}%
Let $f_1,...,f_m$ and $g$ be polynomials in the ring
$F[x_1,...,x_n]$. If each common root of $f_1,...,f_m$ is a root of
$g$, then there exists an integer $r$ such that $g^r$ belongs to the
ideal generated by $f_1,...,f_m$.

\vspace{3mm}\noindent {\bf 89 Riemann-Roch theorem:~~}%
Let $X$ be a curve of genus $g$ and $D$ be a divisor on  $X$. Then
$l(D)-l(K-D)=\deg D +1-g$, where $K$ is the canonical divisor on
$X$.

\vspace{3mm}\noindent {\bf 90 Hurwitz's Theorem:~~}%
Let $f:X\rightarrow Y$ be a finite separable morphism of curves, and
let $n=\deg f$. Let $R$  be the ramification divisor of $f$. Denote
the genus of $X$ and $Y$ by $g(X)$ and $g(Y)$ respectively. Then,
$2g(X)-2=n(2g(Y)-2)+\deg R$.

\vspace{3mm}\noindent {\bf 91 Hasse's theorem:~~}%
If $H$ is the number of points on the elliptic curve $E$ over a
finite field with $q$ elements, then $|H-(q+1)|<2\sqrt{q}$.

\vspace{3mm}\noindent {\bf 92 A theorem of Weil:~~}%
Let  the elliptic curve $E$ be define over a finite field $F_q$ and
$m$ a positive integer. Denote the number of points on the elliptic
curve $E$ over a finite field $F_{q^m}$ by $N$. Then
$N=q^m+1-a^m-b^m$, where $a$ and $b$ satisfy $ab=q$ and
$a+b=H-(q+1)$ with the notation above $H$.

\vspace{3mm}\noindent {\bf 93 R\"{u}ck-Voloch theorem:~~}%
Let the elliptic curve $E$ be define over a finite field $F_q$. Then
the group $E(F_q)$ is isomorphic to a unique direct product of two
cyclic groups $Z_m$ and $Z_n$ with $m|n$ and $m|(q-1)$.

\vspace{3mm}\noindent {\bf 94 A theorem of Mordell:~~}%
For an elliptic curve $E$ over the rational number field $Q$, the
group $E(Q)$ of rational points of $E$ is a finitely-generated
abelian group.

\vspace{3mm}\noindent {\bf 95 Mordell -Weil theorem:~~}%
For an abelian variety $A$ over a number field $K$, the group $A(K)$
of $K$-rational points of $A$ is a finitely-generated abelian group.

\vspace{3mm}\noindent {\bf 96 Faltings' theorem:~~}%
Let $C$ be a non-singular algebraic curve over the rational number
field of genus $g>1$. Then the number of rational points on $C$ is
finite.

\vspace{3mm}\noindent {\bf 97 Tunnell's theorem:~~}%
Let $n$ be a congruent number, if  $n$ is odd then $2A_n=B_n$ and if
 $n$ is even then $2C_n=D_n$, where
 $A_n=|\{(x,y,z)\in Z^3:n=2x^2+y^2+32z^2\}|$,
 $B_n=|\{(x,y,z)\in Z^3:n=2x^2+y^2+8z^2\}|$, $C_n=|\{(x,y,z)\in Z^3:n=8x^2+2y^2+64z^2\}|$, $D_n=|\{(x,y,z)\in
 Z^3:n=8x^2+2y^2+16z^2\}|$.

\vspace{3mm}\noindent {\bf 98 Mazur's torsion theorem:~~}%
The torsion subgroups of the group of rational points on an elliptic
curve defined over the rational number field is one of the following
fifteen groups: $Z/NZ(1\leq N\leq10)$ or $Z/12Z$; $Z/2Z\times Z/2NZ
(1\leq N\leq4)$.

\vspace{3mm}\noindent {\bf 99 Fermat's last theorem:~~}%
If $n>2$ is a positive integer, then the equation $x^n+y^n=z^n$ is
unsolvable in nature number.

\vspace{3mm}\noindent {\bf 100 The Modularity theorem:~~}%
All rational elliptic curves arise from modular forms.

%%%%%%%%%%%%%%%%%%%%%%%%%%%%%%%%%%%%%%%%%%%%%%%%%%%%%%%%%%%%%%%%
%  文章结束
%%%%%%%%%%%%%%%%%%%%%%%%%%%%%%%%%%%%%%%%%%%%%%%%%%%%%%%%%%%%%%%%
\clearpage
\end{document}